\documentclass[11pt]{article}
 \usepackage{benstyle}
 \usepackage{bm}
 
\usepackage[font=small]{caption}
\usepackage{url}

\title{
Frames and numerical approximation
}
\author{Ben Adcock\footnote{Department of Mathematics,	Simon Fraser University, 8888 University Drive, Burnaby, BC V5A 1S6, Canada (\texttt{ben\_adcock@sfu.ca}, \texttt{http://www.benadcock.ca})}
\and 
Daan Huybrechs\footnote{Department of Computer Science, KU Leuven, Celestijnenlaan 200A, BE-3001 Leuven, Belgium (\texttt{daan.huybrechs@cs.kuleuven.be}, \texttt{http://people.cs.kuleuven.be/\textasciitilde daan.huybrechs/})}
}

\begin{document}
\maketitle

\begin{abstract}
Functions of one or more variables are usually approximated with a basis: a complete, linearly-independent system of functions that spans a suitable function space. The topic of this paper is the numerical approximation of functions using the more general notion of frames: that is, complete systems that are generally redundant but provide infinite representations with bounded coefficients.  While frames are well-known in image and signal processing, coding theory and other areas of applied mathematics, their use in numerical analysis is far less widespread.  Yet, as we show via a series of examples, frames are more flexible than bases, and can be constructed easily in a range of problems where finding orthonormal bases with desirable properties (rapid convergence, high resolution power, etc.) is difficult or impossible.  
%These examples highlight three general and useful ways in which frames naturally arise: restrictions of orthonormal bases to subdomains, augmentation of orthonormal bases by a finite number of `feature' functions, and concatenation of two or more orthonormal bases.  
For instance, we exhibit a frame which yields simple, high-order approximations of smooth, multivariate functions in arbitrary geometries.

A key concern when using frames is that computing a best approximation requires solving an ill-conditioned linear system.  Nonetheless, we construct a frame approximation via regularization with bounded condition number (with respect to perturbations in the data), and which approximates any function up to an error of order $\sqrt{\epsilon}$, or even of order $\epsilon$ with suitable modifications.  Here $\epsilon$ is a threshold value that can be chosen by the user.Ê Crucially, rate of decay of the error down to this level is determined by the existence of approximate representations of $f$ in the frame possessing small-norm coefficients.  We demonstrate the existence of such representations in all of our examples.  Overall, our analysis suggests that frames are a natural generalization of bases in which to develop numerical approximation. In particular, even in the presence of severely ill-conditioned linear systems, the frame condition imposes sufficient mathematical structure in order to give rise to accurate, well-conditioned approximations.
\end{abstract}

\paragraph{Keywords}frames, function approximation, ill-conditioning, singular value decomposition

\paragraph{AMS subject classifications} 42C15, 42C30, 41A10, 65T40

\section{Introduction}
Frames are a generalization of bases that allow for redundancy amongst the generating elements.  They are indispensable tools in modern signal and image processing,  and are widely used in a range of other problems, such as compression, source coding, robust transmission and sampling theory \cite{BenedettoIrregular,christensen2003introduction,KovacevicChebiraFrames2,KovacevicChebiraFrames1,mallat09wavelet}.  Yet frames -- specifically, frames of functions --  are generally less well known in numerical analysis.\footnote{Finite frames of vectors, on the other hand, are better known in  the numerical linear algebra community (see, for instance, \cite{DattaEtAlETF,FickusEtAlETF,RenesETF,StrohmerETF,SustikEtAlETF} and references therein).  However, the focus of this article is the infinite setting of frames of function spaces.} Although they arise in approximation problems in a number of ways, a systematic and general study of numerical frame approximation does not appear to have been undertaken.

The purpose of this paper is to consider frames from this perspective.  By means of motivation, we introduce three classes of problems in numerical computing where frames already occur naturally, or where they may potentially lead to better methods.  Our main objective is to examine the accuracy and conditioning of frame approximations, and the properties of a given frame which affect both.

A key theme of this paper is the difference between the behaviour of infinite frames and the corresponding truncated frames used in approximation.  Quite unlike the case of orthonormal bases, the matrices of the linear systems associated with best approximations in frames (truncated Gram matrices) have condition numbers that necessarily grow with the degree of the approximation.  Moreover, this growth can be arbitrarily fast.  The surprising fact, however, is that accurate frame approximations can still be computed, for instance via suitable regularization.
Crucially, and unlike in the case of orthonormal bases, this means that there are fundamental differences between `theoretical' frame approximations (e.g.\ the best approximation) and `numerical' frame approximations (i.e.\ the solution of the regularized problem).  Understanding and documenting these differences is a central aspect of this paper.

\subsection{Orthonormal bases}
Let $\Phi  = \{ \phi_n \}_{n \in \bbN}$ be an orthonormal basis of a Hilbert space $\rH$, indexed for convenience over the natural numbers $\bbN = \{1,2,\ldots \}$. 

Two key properties of $\Phi$ are as follows.  First, the straightforward representation of any $f \in \rH$ in the basis using the inner product on $H$.  That is,
\be{
\label{ONBexpansion}
f = \sum_{n \in \bbN} \ip{f}{\phi_n} \phi_n,\qquad \forall f \in \rH,
}
where the infinite sum converges in $\rH$.  Note that the \textit{coefficients} $\{ \ip{f}{\phi_n} \}_{n \in \bbN}$ in this expansion are unique.  Second, Parseval's identity, which asserts that the Hilbert space norm of $f$ is precisely the $\ell^2$-norm of its coefficients:
\be{
\label{ParsevalONB}
\| f \|^2 = \sum_{n \in \bbN} | \ip{f}{\phi_n} |^2,\qquad \forall f \in \rH.
}
If the coefficients $\{ \ip{f}{\phi_n} \}_{n \in \bbN}$ are known (or have been computed), approximation in $\Phi$ is a straightforward affair.  One simply replaces \R{ONBexpansion} by a finite expansion
\be{
\label{ONBorthproj}
f \approx \sum_{n=1}^{N} \ip{f}{\phi_n} \phi_n.
}
This approximation has the beneficial property of being the orthogonal projection onto the finite-dimensional subspace $\rH_N = \spn \{ \phi_n \}_{n=1}^{N}$, and therefore the best approximation to $f$ from $\rH_N$ in the norm of $\rH$ (we assume throughout that $\rH_N$ is endowed with this norm).

\subsection{Frames}
An indexed family $\Phi  = \{ \phi_n \}_{n \in \bbN}$ is called a \textit{frame} for $\rH$ if it satisfies the so-called \emph{frame condition}
\be{
\label{framecondintro}
A \| f \|^2 \leq \sum_{n \in \bbN} | \ip{f}{\phi_n} |^2 \leq B \| f \|^2,\quad \forall f \in \rH,
}
for constants $A,B > 0$.  The optimal constants $A,B>0$ such that \R{framecondintro} holds, i.e.\ the largest possible $A$ and the smallest possible $B$, are referred to as the \textit{frame bounds}~\cite{christensen2003introduction,duffin1952frames}.    We recall frames and their theory in more detail in \S\ref{sect:preliminaries}.

Generalizing Parseval's identity \R{ParsevalONB}, the frame condition expresses a norm equivalence between the $\ell^2$-norm of the coefficients $\{ \ip{f}{\phi_n}\}_{n \in \bbN}$ and the Hilbert space norm of $f$.  Yet frames differ from orthonormal bases in a number of key ways:
\begin{itemize}
\item[(i)] The frame elements $\phi_n$ are not generally orthogonal.
\item[(ii)] While \R{framecondintro} implies that $\spn(\Phi)$ is dense in $\rH$, $\Phi$ need not be a basis.  In fact, a frame is typically \textit{redundant}.  That is, for any $f \in \rH$ it is possible to find more than one sequence of square-summable coefficients $\bm{c} = \{ c_n \}_{n \in \bbN}$ such that $f = \sum_{n \in \bbN} c_n \phi_n$.
\item[(iii)] In general, a representation such as \R{ONBexpansion} does not hold.
\end{itemize}

\subsection{Computing orthogonal projections with frames}
While (i) means that frames are more flexible than orthonormal bases (indeed, \R{framecondintro} is far less restrictive a condition than orthogonality), it presents an immediate difficulty for numerical approximation with frames.  Even if the coefficients $\ip{f}{\phi_n}$ are known, the orthogonal projection cannot be expressed as in \R{ONBorthproj} as a sum of these coefficients multiplied by the frame elements.  Determining the coefficients $\bm{x} = \{ x_n \}^{N}_{n=1}$ of the orthogonal projection with respect to the frame elements requires solving a particular linear system (see \S \ref{ss:BestApprox} for further details).  In the notation of our paper, this is denoted by
\be{
\label{first_system}
\bm{G}_N \bm{x} = \bm{y},\qquad \bm{y} = \{ \ip{f}{\phi_n} \}^{N}_{n =1},
}
where $\bm{G}_N$ is the $N \times N$ \emph{truncated Gram matrix}
\[
\bm{G}_N = \left \{ \ip{\phi_m}{\phi_n} \right \}^{N}_{n,m =1} \in \bbC^{N \times N}.
\]
If $\bm{x} \in \bbC^{N}$ is a solution of \R{first_system}, the orthogonal projection is given by $\sum^{N}_{n=1} x_n \phi_n$.

The second observation (ii) raises another, and arguably more crucial, issue in practical frame approximation.  Due to orthogonality, the Gram matrix of a finite subset $\{ \phi_n \}^{N}_{n=1}$ of an orthonormal basis is perfectly conditioned; indeed, $\bm{G}_N = \bm{I}$ is the identity matrix.  It is tempting to think that the frame condition \R{framecondintro} endows the Gram matrices of finite subsets of frames with a similar property.  However, while the \textit{truncated} frame $\Phi_N = \{ \phi_n \}^{N}_{n=1}$ is indeed a frame for its span $\rH_N$, and thus satisfies a frame condition, its frame bounds $A_N$ and $B_N$ may behave wildly as $N \rightarrow \infty$, even when the infinite frame bounds $A$ and $B$ are mild.  We shall see examples later in this paper where the ratio $B_N/A_N$ grows exponentially fast in $N$.  If $\Phi_N$ is linearly independent (a condition which is satisfied by all examples of this paper -- see \S \ref{ss:frames}), the condition number of the Gram matrix $\bm{G}_N$ is precisely $B_N/A_N$ (Lemma \ref{l:GramNormFrameBds}).  Hence ill-conditioning is equivalent to poorly-behaved frame bounds.  Understanding this ill-conditioning and its effect on the resulting numerical frame approximation obtained by solving a regularized version of \R{first_system} is the central theme of this paper.

We note in passing that this ill-conditioning stems from the noninvertibility of the Gram \textit{operator} $\cG$ of the frame $\Phi$ (see \S \ref{ss:frames}). Such noninvertibility is due to (ii).  Discretizing $\cG$ with the matrix $\bm{G}_N$ results in small, nonzero 
eigenvalues approximating the zero eigenvalue of $\cG$ \cite{HarrisonThesis}, and therefore large condition numbers (see Remark \ref{r:opthy}).

\subsection{Motivations}

Orthonormal bases are ubiquitous in numerical analysis.  Important cases include Fourier and Chebyshev bases, in which case fast algorithms exist to (approximately) compute the expansions \cite{trefethen2013atap}.  A major disadvantage of orthogonal bases however is their inflexibility.  As an example, consider the problem of approximating smooth functions of one or more variables.  While it is easy to construct good\footnote{The word `good' in this paper is taken to mean \textit{spectrally} convergent, i.e.\ having rates of convergence depending only on the smoothness of the function being approximated.} orthogonal bases of functions on intervals,  it is much harder to do so in higher dimensions unless the domain is particularly simple (e.g.\ a hypercube). This aside, it is also problematic to find a good basis for singular functions, or to force periodicity on nonperiodic problems in order to take advantage of the FFT.

In this paper we show that good frames can be found for all these problems.  In particular, we identify a simple frame with spectral rates of convergence (in the corresponding Hilbert space norm) for approximating functions defined on arbitrary Lipschitz domains; see Fig.\ \ref{f:northpole} below. These examples illustrate three different generic constructions which always lead to frames:  restrictions of orthonormal bases to subdomains, augmentation of an orthonormal basis by a finite number of additional terms, and concatenation of several orthonormal bases. This leads us to opine that frames are natural tools for many problems in numerical analysis where constructing orthonormal bases is difficult or impossible.

\subsection{Overview and main results}

We restrict our focus in this paper to two key properties: the convergence of frame approximations and their conditioning.  We shall mostly ignore the question of efficiency, since this is highly dependent on the type of frame used and in this paper we strive for generality (we return to this topic briefly in \S \ref{s:conclusions}).
Our main conclusion is the following.  In spite of the extreme ill-conditioning of the linear system \R{first_system}, accurate frame approximations, in a sense we make precise below, can be computed numerically.  To do so, the linear system \R{first_system} is regularized using a truncated Singular Value Decomposition (SVD) of the Gram matrix $\bm{G}_N$.  

We recall the main elements of frame theory in \S\ref{sect:preliminaries}. In \S\ref{sect:examples} we introduce three generic constructions of frames that are useful in numerical approximations along with our three main examples.  All these examples deal with the problem of approximating functions where it is not straightforward or even desirable to use orthonormal bases, and where frame approximations present a viable alternative.

Our analysis commences in \S\ref{sect:IllCond}. The ratio $B/A$ of the frame bounds is an important quantity in the conditioning of frames.  However, unlike for orthogonal or Riesz bases, passing from the countable frame $\Phi$ to a finite subset $\Phi_N = \{ \phi_n \}^{N}_{n=1}$ necessarily causes a deterioration in the frame bounds.  We document this phenomenon in \S\ref{sect:IllCond}.  Lemma~\ref{lem:pollution} shows how frame bounds deteriorate after truncation, and Proposition~\ref{p:simple_bad_frame} establishes that the effect can be arbitrarily bad. The condition numbers $\kappa(\bm{G}_N)$ for the three example frames are estimated in \S\ref{ss:example_condition_numbers}. 

%Each exhibits algebraic, superalgebraic or exponential growth in $N$.

We consider the computation of the best approximation via orthogonal projection in \S\ref{sect:OP_compute}. We first show in Proposition~\ref{p:unbounded_solution} that the $\ell^2$-norm of the \emph{exact} solution $\bm{x}$ of the system \R{first_system} is generally unbounded in $N$, due to the ill-conditioning. Hence, computing the best approximation with any accuracy for large $N$ in floating point arithmetic is typically impossible. However, the situation improves markedly after regularizing $\bm{G}_N$ by truncating its singular values below a threshold $\epsilon$. In Theorem~\ref{t:Proj_err} we show that the convergence of the resulting regularized projection to a function $f$ is dictated by how well $f$ can be approximated by vectors of coefficients with small norm. Specifically, if $\cP^{\epsilon}_N f$ is the regularized projection, then
\be{
\label{fn_err}
\| f - \cP^{\epsilon}_N f \| \leq \inf \left \{ \nm{ f - \cT_N \bm{z}} + \sqrt{\epsilon} \| \bm{z} \| : \bm{z} \in \bbC^N \right \},
}
where $\cT_N \bm{z} =  \sum^{N}_{n=1} z_n \phi_n \in \rH_N$.  The first term $\nm{f - \cT_N f}$ is standard and represents the approximation error corresponding to a coefficient vector $\bm{z}$.  The second term $\sqrt{\epsilon} \nm{\bm{z}}$ is uncommon in the literature on frames, and indeed it is specific to a numerical frame approximation.  Its effect is to limit the overall approximation error $\nm{ f - \cT_N \bm{z}} + \sqrt{\epsilon} \| \bm{z} \|$ by penalizing larger-norm coefficient vectors $\bm{z}$. In other words, the overall approximation error will only be small asymptotically if $f$ can be represented in the frame (first term) with coefficients of small norm (second term). The existence of such representations in the first place is guaranteed by the frame property, and that is why we argue that the mathematical structure of a frame seems a highly appropriate general context in which to discuss function approximation in redundant systems.

Theorem~\ref{t:coefficients1} shows that the coefficient vector $\bm{x}^{\epsilon}$ of the regularized projection $\cP^{\epsilon}_N f$ eventually (for large $N$) exhibits small norm, though there may be an initial regime in which it is large. The precise result is as follows:
\be{
\label{coeff_err}
\| \bm{x}^{\epsilon} \| \leq \inf \left \{ 1/\sqrt{\epsilon} \| f - \cT_N \bm{z} \| + \| \bm{z} \| : \bm{z} \in \bbC^N \right \}.
}
Again, the frame property is crucial here, in that it ensures the existence of vectors of coefficients $\bm{z}$ with small norm such that $\cT_N \bm{z} \approx f$.  The price to pay for this beneficial property is that the true convergence rate of the best approximation may not be realized after regularization. Instead, one finds best approximations subject to having a small-norm coefficient vector. In practice, however, small-norm coefficient vectors are often more desirable.

A point of clarification. At this stage, the reader may be tempted to conclude that frame approximations are of limited use in practice, since they can obtain at best $\ord{\sqrt{\epsilon}}$ accuracy.  Furthermore, Theorem \ref{t:Proj_Condit} illustrates that the mapping $\bm{y} \mapsto \cP^{\epsilon}_N f$ from the data $\bm{y}$ (see \R{first_system}) to the approximation $\cP^{\epsilon}_N f$ is not completely well-conditioned, since its (absolute) condition number can be behave like $\cO(1/\sqrt{\epsilon})$.  However, these issues can be simultaneously overcome.  In \S \ref{s:numstableapprox} we briefly describe a generalized frame approximation which achieves $\ord{\epsilon}$ accuracy, and has a bounded condition number, independent of $\epsilon$.  The full analysis of these techniques (which builds on this paper) is described in a companion paper \cite{BADHFramesPart2}.

\begin{figure}[t]
\begin{center}
 \includegraphics[width=10cm]{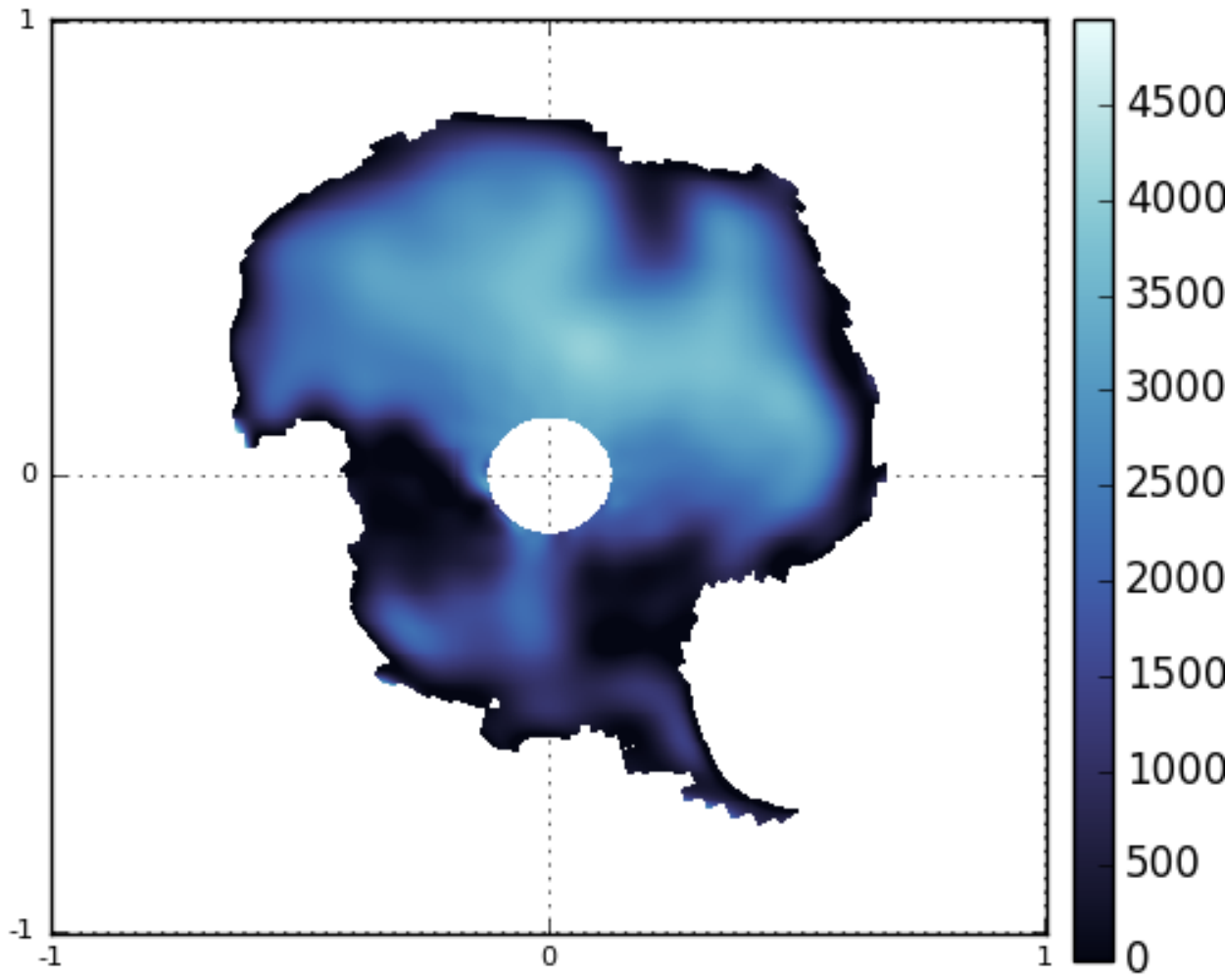}
\end{center}
\caption{Illustration of the elevation of ice at the north pole (in $meters$). Discrete height data from EOS, NASA's Earth Observing System, was approximated in a discrete least-squares sense using using a conventional two-dimensional Fourier series on a bounding box (the data is publicly available on NASA's EOS website at \texttt{https://neo.sci.gsfc.nasa.gov/}). The data was supplied in equispaced points in polar coordinates, and the approximation was constructed in polar coordinates as well. No data was supplied outside the north pole, i.e.\ the function samples implicitly define the shape of the domain. Data was also absent beyond a certain latitude, hence the circle of missing data near the pole. The irregular and punctured shape of the domain, without explicitly defined boundary, does not hinder the numerical frame approximation. The truncation parameter $\epsilon$ was chosen to be fairly large, in order to reduce the growth of the coefficients while maintaining proximity to the data, as suggested by Theorems \ref{t:Proj_err} and \ref{t:coefficients1}. The algorithm used is described and analyzed in \cite{matthysen2017fastfe2d}. It exploits the equispaced property of the data points via the FFT, even though data is only given on an irregular subset of the full rectangular grid. In this paper, we focus not on the algorithm, but on the mathematical analysis that make this type of approximation possible and reliable.}
\label{f:northpole}
\end{figure}

Finally, let us briefly demonstrate the effectiveness of numerical frame approximations.  Fig.\ \ref{f:northpole} shows the approximation of ice elevation at the north pole from discrete height data.  The domain is highly irregular, and data is also missing at high latitudes.  Yet an accurate approximation is achieved using a simple and spectrally-accurate frame, whose elements $\phi_n$ are nothing more than the usual Fourier basis functions on the bounding box $[-1,1]^2$ (Example 1 of this paper).

\subsection{Relation to existing work}\label{ss:background}

Frames were introduced in the context of nonharmonic Fourier series by Duffin and Schaeffer \cite{duffin1952frames}.  They were later developed in the 1980s by Daubechies, Grossmann and Meyer \cite{DaubechiesGrossmannMeyerPainless} with the systematic study of wavelets.  Since then they have become an integral part of modern signal and image processing, compression, coding theory and sampling theory.  For overviews, see \cite{BenedettoIrregular,CasazzaArtFrameThy,FiniteFramesBook,christensen2003introduction,daubechies1992tenlectures,KovacevicChebiraFrames2,KovacevicChebiraFrames1}.

With the notable exceptions of \cite{HarrisonThesis,StrohmerNonuniform}, \cite{SongEtAlFrame,SongGelbInverse} and \cite{Dahlke2007,Stevenson2003}, approaches to numerical frame approximations have usually centered around either explicitly identifying an appropriate \textit{dual} frame or numerically inverting the frame \textit{operator} (from which the so-called \textit{canonical dual} frame can be computed) \cite{CasazzaChristensenFrameCoeffApprox,ChristensenCasazzaJAA,CasazzaChristensenJAT,Christensen1993Projection,Christensen1996FramesRieszBasis,ChristensenJFAA,ChristensenLindnerFrameRB,ChristensenStrohmerFS}.  We refer to \S \ref{ss:canonical_dual_frame} for further details. While these approaches are useful for approximations with, for example, Gabor or wavelet frames and their various generalizations (e.g.\ multi-wavelets \cite{christensen2003introduction}, ridgelets \cite{emmanuel1999ridgelets}, curvelets \cite{candes2004new,mallat09wavelet} and shearlets \cite{shearletsbook}), for the problems which motivate this paper the dual frame expansion usually converges too 
slowly to be of practical use.  We give several examples of this phenomenon later.  Conversely, our focus in this paper is on computing best approximations with frames, or more precisely, surrogates obtained from solving regularized systems. The regularized Gram systems we consider in this paper have previously been studied in \cite{HarrisonThesis,StrohmerNonuniform} in the context of frames of exponentials arising in nonuniform sampling problems.  In particular, \cite[Thm.\ 5.17]{HarrisonThesis} asserts convergence of the coefficients $\bm{x}^{\epsilon}$ to the so-called frame coefficients (see \S \ref{ss:canonical_dual_frame}) as $N \rightarrow \infty$.  We extend this work by establishing the convergence rate \R{fn_err} and coefficient bound \R{coeff_err} for arbitrary frames.  Note that our regularization procedure is also related to a recent approach of \cite{DorflerGabor}, used for the construction of adaptive Gabor frame approximations.

Our study of numerical frame approximation stems from previous works on so-called \textit{Fourier extensions} \cite{BoydFourCont,brunoFEP}, also known as \textit{Fourier continuation} or \textit{Fourier embedding} in the context of numerical PDEs \cite{pasquettiFourEmbed}.  The connection to frame theory was first explored by the second author in \cite{huybrechs2010fourier}, and further developed in the one-dimensional setting by both authors in \cite{FEStability}.  A by-product of this paper is an extension of the results of \cite{FEStability} to $d \geq 1$ dimensions.  Yet we stress that our main results apply to any frame, not just Fourier extensions.

We also draw several interesting connections to other fields.  In particular, the disparity between truncated frames and infinite frames is related to the spectral theory of self-adjoint operators, and specifically the phenomenon of pollution in the finite section method \cite{DaviesPlumSpecPoll,LewinSereSpecPoll}.  We also make links to the topic of time- and band-limiting \cite{HoganLakeyBook}, in particular the prolate spheriodal wavefunctions \cite{SlepianV}, and to classical regularization theory of ill-posed problems \cite{EnglRegularization,HansenEtAlLeastSquares,NeumaierIllCond}.

\subsection{On supplementary materials}
This paper is accompanied by supplementary materials \cite{framespart1SM}.  These materials contain the proofs of a number of results pertaining to the various frames used as examples in the paper.  Conversely, all results that are general (that is, apply to any frame) are proved in the main manuscript.

\section{Preliminaries}
\label{sect:preliminaries}

\subsection{Orthogonal and Riesz bases}
For the remainder of the paper, $\Phi = \{ \phi_n \}_{n \in I}$ is an indexed family in a separable Hilbert space $\rH$ over the field $\bbC$, where $I$ is a countable index set (for instance, the natural numbers $\bbN$ or integers $\bbZ$).  We write $\ip{\cdot}{\cdot}$ and $\nm{\cdot}$ for the inner product and norm on $\rH$ respectively.  

The indexed family $\Phi$ is an \textit{orthonormal} basis for $\rH$ if $\spn(\Phi)$ (the vector space of all finite linear combinations of elements of $\Phi$) is dense in $\rH$ and $\ip{\phi_n}{\phi_m}$ = $\delta_{n,m}$,$\forall n,m \in I$.  Recall that orthonormal bases satisfy Parseval's identity
\be{
\label{Parseval2}
\| f \|^2 = \sum_{n \in I} | \ip{f}{\phi_n} |^2,\quad \forall f \in \rH. 
}
Equivalently,
\be{
\label{Parseval1}
\| \bm{x} \| = \nmu{\sum_{n \in I} x_n \phi_n },\quad \forall \bm{x} = \{ x_n \}_{n \in I} \in \ell^2(I).
}
Here and throughout, $\ell^2(I)$ denotes the space of square-summable sequences indexed over $I$, and $\nm{\cdot}$ denotes its norm, i.e.\ $\nm{\bm{x}} = \sqrt{\sum_{n \in I} |x_n|^2 }$.  

\textit{Riesz} bases are generalizations of orthonormal bases.  $\Phi$ is a Riesz basis of $\rH$ if $\spn(\Phi)$ is dense in $\rH$ and the following relaxed version of \R{Parseval1} holds:
\be{
\label{RieszParseval1}
A \| \bm{x} \|^2 \leq \nmu{\sum_{n \in I} x_n \phi_n }^2 \leq B \| \bm{x} \|^2,\quad \forall \bm{x} = \{ x_n \}_{n \in I} \in \ell^2(I).
}
Here $A,B > 0$ are positive constants.  Throughout this paper, whenever constants $A$ and $B$ are introduced in an inequality such as this, they will be taken to be the optimal constants such that the inequality holds.  Note that this inequality also implies the following relaxed version of \R{Parseval2}:
\be{
\label{RieszParseval2}
A \| f \|^2 \leq  \sum_{n \in I} | \ip{f}{\phi_n} |^2 \leq B \| f \|^2,\quad \forall f \in \rH.
}
Every orthonormal basis is a Riesz basis, but Riesz bases need not be orthonormal, or even orthogonal.  A popular example of such a basis are the hat functions commonly used in finite element methods. They are not orthogonal, but they are a basis for their span. Another example are more general B-splines \cite{christensen2003introduction}.

Any Riesz basis has a unique \textit{dual} Riesz basis $\Psi$.  This basis satisfies
\bes{
\ip{\phi_n}{\psi_m} = \delta_{n,m}, \quad n,m \in I.
}
For this reason, the Riesz basis and its dual are termed \emph{biorthogonal}. From this, the unique representation of an element $f \in \rH$ in the Riesz basis is given explicitly in terms of inner products with the dual basis.  That is,
\be{
 f = \sum_{n \in I} \ip{f}{\psi_n} \phi_n = \sum_{n \in I} \ip{f}{\phi_n} \psi_n,
}
with convergence in $\rH$. Note that an orthonormal basis is self-dual, i.e. $\Psi = \Phi$.

\subsection{Frames}\label{ss:frames}
$\Phi$ is called a \textit{frame} for $\rH$ if it satisfies the \emph{frame condition}, which we repeat here for convenience:
\be{
\label{framecond}
A \| f \|^2 \leq \sum_{n \in I} | \ip{f}{\phi_n} |^2 \leq B \| f \|^2,\quad \forall f \in \rH.
}
Note that \R{framecond} implies that $\spn(\Phi)$ is dense in $\rH$.  It follows from the Parseval identity \R{Parseval2} and its generalization \R{RieszParseval2} that orthonormal and Riesz bases are frames.  However, most frames are not bases.  Indeed, frames are generally not \textit{$\omega$-independent} (see  \cite[Sec.\ 6.1]{christensen2003introduction}, for example): that is, there exist nonzero coefficients $\{ x_n \}_{n \in I}$ for which the sum $\sum_{n \in I} x_n \phi_n$ converges in $\rH$ and satisfies $\sum_{n \in I} x_n \phi_n = 0$.  Conversely, bases are always $\omega$-independent.
As mentioned, this \textit{redundancy} gives frames far greater flexibility than bases, making them easier to construct for particular problems.

We now introduce several standard concepts in frame theory.  Associated to any frame $\Phi$ is the so-called \textit{synthesis} operator
\bes{
\cT : \ell^2(I) \rightarrow \rH,\quad \bm{y} = \{ y_n \}_{n \in I} \mapsto \sum_{n \in I} y_n \phi_n.
}
Its adjoint, the \textit{analysis} operator, is given by
\bes{
\cT^* : \rH \rightarrow \ell^2(I),\quad f \mapsto \{ \ip{f}{\phi_n} \}_{n \in I},
}
and the composition $\cS = \cT \cT^*$, known as the \textit{frame} operator, is 
\bes{
\cS : \rH \rightarrow \rH,\quad f \mapsto \sum_{n \in I} \ip{f}{\phi_n} \phi_n.
}
Note $\cS f$ converges to $f$ in $\rH$ for an orthonormal basis, but for a Riesz basis or frame this is no longer the case in general.  Still, the frame operator is a useful object. It is self-adjoint by construction, and it follows from the frame condition that $\cS$ is also bounded and invertible on $\rH$ \cite[Lemma 5.1.5]{christensen2003introduction}.

The Gram operator of a frame is defined by $\cG = \cT^* \cT$.  That is,
\be{
\label{GramOp}
\cG : \ell^2(I) \rightarrow \ell^2(I),\ \bm{x} = \{x_n \}_{n \in I} \mapsto \left \{ \sum_{m \in I} \ip{\phi_m}{\phi_n} x_m \right \}_{n \in I}.
}
Note that $\cG$ is a bounded operator on $\ell^2(I)$, but is not in general invertible (see \S \ref{sect:IllCond}).  We may also view $\cG$ as the infinite matrix $\bm{G} = \left \{ \ip{\phi_n}{\phi_m} \right \}_{n,m \in I}$.  Throughout this paper all infinite matrices are equivalent to bounded operators on $\ell^2(I)$.

A frame is said to be \textit{tight} if $A=B$, in which case $\cS = A \cI$ is a multiple of the identity $\cI$.  However, the corresponding Gram operator $\cG$ of a tight frame is not a multiple of the identity, unless the frame happens to be an orthonormal basis (for its closed linear span).

We shall also need two further notions.  First, a frame is said to be \textit{exact} if it ceases to be a frame when any one element is removed.  A frame that is not exact is referred to as \textit{inexact}.  Second, we say a frame $\{ \phi_n \}_{n \in I}$ is \textit{linearly independent} if every finite subset $\{ \phi_n \}_{n \in J}$, $|J| < \infty$, is linearly independent. 

A frame is exact if and only if it is a Riesz basis \cite[Theorem 5.5.4]{christensen2003introduction}.  Hence, for the remainder of this paper we will assume that all frames are inexact.  We shall also assume that all frames are linearly independent.  This is mainly for convenience, and it will be the case in all examples discussed.  Note that a linearly-independent frame is not necessarily a Riesz basis.  See \cite[Chpt.\ 6]{christensen2003introduction} for further discussion on independence and the relations between frames and Riesz bases.

\subsection{Dual frames}
\label{ss:canonical_dual_frame}

A frame $\Psi = \{ \psi_n \}_{n \in I} \subseteq \rH$ is called a \textit{dual} frame for $\Phi$ if
\be{
\label{dual_frame}
f = \sum_{n \in I} \ip{f}{\psi_n} \phi_n = \sum_{n \in I} \ip{f}{\phi_n} \psi_n,\qquad \forall f \in \rH.
}
An inexact frame necessarily has more than one dual frame.  Moreover, a frame and its duals are not biorthogonal, unlike in the case of Riesz bases.
However, there is a unique so-called \textit{canonical} dual frame $\Psi = \{ \psi_n \}_{n \in I}$, with elements given by $\psi_n = \cS^{-1} \phi_n$ where $\cS$ is the frame operator. Since $\Psi$ is a dual frame, one has
\be{
\label{dual_rep}
f = \sum_{n \in I} \ip{f}{\cS^{-1} \phi_n} \phi_n = \sum_{n \in I} \ip{\cS^{-1} f}{\phi_n} \phi_n.
}
Furthermore, the canonical dual frame bounds are $1/B$ and $1/A$ respectively, i.e.\
\be{
\label{dual_bd}
1/B \| f \|^2 \leq \sum_{n \in I} | \ip{f}{\cS^{-1} \phi_n} |^2 \leq 1/A \| f \|^2,\quad \forall f \in \rH.
}
We refer to the coefficients $\bm{a} = \{ \ip{f}{\cS^{-1} \phi_n } \}_{n \in I}$ as the \textit{frame coefficients} of $f$.  
Note that these coefficients have the beneficial property that, amongst all possible representations of $f$ in $\Phi$, they have the smallest norm. That is, if $f = \sum_{n \in I} a_n \phi_n = \sum_{n \in I} c_n \phi_n$ for some $\bm{c} = \{ c_n \}_{n \in I}$, then $\| \bm{c} \| \geq \| \bm{a} \|$ \cite[Lem.\ 5.4.2]{christensen2003introduction}.

At this stage, one might be tempted to approximate $f$ by computing its dual frame coefficients and truncating the expansion \R{dual_frame}.  As noted in \S \ref{ss:background}, this could potentially be done either by analytically identifying a dual frame (when possible) or by numerically inverting the frame operator.  However, computational issues aside -- in the case where the frame is not tight, this requires inversion of the operator $\cS$, which has infinite-dimensional domain and range -- the approximation $\sum_{n \in I_N} \ip{\cS^{-1} f}{\phi_n} \phi_n$ is generally not the orthogonal projection onto $\rH_N = \spn \{ \phi_n : n \in I_N \}$.  For the examples which motivate this paper, this approximation typically converges much more slowly in the norm of $\rH$ than the orthogonal projection.\footnote{This is in contrast to the case of wavelet frames and their various generalizations, which are specifically designed to have accurate dual frame representations \cite{christensen2003introduction,mallat09wavelet}.}  See \S \ref{sect:examples} and Figs.\ \ref{f:FE_Err_Coeff} and \ref{f:FP_Err_Coeff} for several examples of this phenomenon.

\subsection{Truncated frames}

For each $N \in \bbN$ we introduce the truncated system $\Phi_N = \{ \phi_n \}_{n \in I_N}$ where $I_N \subseteq I$ is a finite index set with $|I_N|= N$.  For simplicity, we assume that the index sets $\{ I_N \}_{N \in \bbN}$ are nested and satisfy
\be{
\label{indexsets}
I_1 \subseteq I_2 \subseteq \ldots, \qquad \bigcup^{\infty}_{N=1} I_N = I.
}
The system $\Phi_N$ is a frame for its span $\rH_N = \mathrm{span}(\Phi_N)$.  We write $A_N,B_N > 0$ for the frame bounds, so that
\be{
\label{TruncFrameBounds}
A_N \| f \|^2 \leq \sum_{n \in I_N} | \ip{f}{\phi_n} |^2 \leq B_N \| f \|^2,\quad \forall f \in \rH_N,
}
and let 
\ea{
\cT_N &: \bbC^N \rightarrow \rH_N, \ \bm{y} = \{ y_n \}_{n \in I_N} \mapsto \sum_{n \in I_N} y_n \phi_n \nn,
\\
\cT^*_N &: \rH_N \rightarrow \bbC^N,\ f \mapsto \{ \ip{f}{\phi_n} \}_{n \in I_N}, \label{TN_def}
\\
\cS_N = \cT_N \cT^*_N &: \rH_N \rightarrow \rH_N,\ f \mapsto \sum_{n \in I_N} \ip{f}{\phi_n} \phi_n, \nn
}
be the truncated synthesis, analysis and frame operators respectively.  We also define the truncated Gram operator $\cG_N = \cT^*_N \cT_N$ and the associated $N \times N$ Gram matrix
\be{
\label{GramMat}
\bm{G}_N = \left \{ \ip{\phi_m}{\phi_n} \right \}_{n,m \in I_N} \in \bbC^{N \times N}.
}
Since the frame $\Phi$ is linearly independent by assumption, it follows that the Gram matrix $\bm{G}_N$ is nonsingular.  Indeed, if $\bm{x} = \{x_n\}_{n \in I_N} \in \bbC^N$ then $\bm{x}^* \bm{G}_N \bm{x} = \nm{ \sum_{n \in I_N} x_n \phi_n }^2$ and,  by linear independence, the right-hand side is zero if and only if $\bm{x} = \bm{0}$.

\subsection{Best approximations and rates of convergence}
\label{ss:BestApprox}

A key task in approximation is to compute the orthogonal projection $\cP_N$ onto $H_N = \mathrm{span}(\Phi_N)$. Observe that $\cP_N f$ is the best approximation to $f$ from $\rH_N$ in $\nm{\cdot}$. For $f \in \rH$, write
\be{
\label{frameOP}
\cP_N f = \sum_{n \in I_N} x_n \phi_n,\quad \bm{x} = \{ x_n \}_{n \in I_N} \in \bbC^N.
}
Since $\cP_N f$ is defined by the orthogonality conditions
\bes{
\ip{\cP_N f}{\phi_n} = \ip{f}{\phi_n}, \quad \forall n \in I_N, 
}
the coefficients $\bm{x} = \{ x_n \}_{n \in I_N}$ are the unique solution of the linear system
\be{
\label{GramLinSyst}
\bm{G}_N \bm{x} = \bm{y},\quad \bm{y} = \{ \ip{f}{\phi_n} \}_{n \in I_N}.
}
Hence, computing the best approximation in a frame requires solving an $N \times N$ linear system.  This system turns out to be ill-conditioned, in direct contrast with the case of an orthonormal basis, wherein the Gram matrix $\bm{G}_N$ is the identity and $x_n = y_n = \ip{f}{\phi_n}$.  The treatment of this ill-conditioning is central task of this paper.

Besides conditioning, another concern of this paper is the convergence rate of the approximation $\cP_N f$ and its surrogates obtained via regularization.  To this end, we distinguish three types of convergence of an approximation $f_N$ to a function $f$.  First, we say that $f_N$ converges \textit{algebraically} fast to $f$ at rate $k$ if $\| f - f_N \| = \ordu{N^{-k}}$ as $N \rightarrow \infty$.  Second, if $\nm{f-f_N}$ decays faster than any algebraic power of $N^{-1}$ then we say that $f_N$ converges \textit{superalgebraically} fast to $f$.  Third, we say that $f_N$ converges \textit{geometrically} fast to $f$ if there exists a $\rho > 1$ such that $\nm{f-f_N} = \ordu{\rho^{-N}}$.

\section{Examples of frames}
\label{sect:examples}

We now introduce three examples of frames that will be used throughout the paper to interpret our main results.  These examples also serve to illustrate the flexibility gained by allowing redundancy in the approximation.

\vspace{1pc} \noindent
\textbf{Example 1.\ Fourier frames for complex geometries.}
Let $\Omega \subseteq \bbR^d$ be a compact domain and $f : \Omega \rightarrow \bbR$.  Besides simple domains (cubes, toruses, spheres, etc.), it is in general difficult to find orthonormal bases for $\rH = \rL^2(\Omega)$ with simple, explicit expressions and whose orthogonal projections converge \textit{spectrally} fast: that is, algebraic for functions with finite orders of smoothness and superalgebraic for infinitely-smooth functions.  However, it is straightforward to find a frame with this property.

Since $\Omega$ is compact, it can be contained in an open hypercube.
Without loss of generality, suppose that $\Omega \subseteq (-1,1)^d$.  Now consider a system of functions formed by the restriction of the orthonormal Fourier basis on $(-1,1)^d$ to $\Omega$:
\be{
\label{FE_frame}
\Phi = \{ \phi_{\bm{n}} \}_{\bm{n} \in \bbZ^d},\qquad \phi_{\bm{n}}(\bm{t}) = 2^{-d/2} \E^{\I \pi \bm{n} \cdot \bm{t}},\quad \bm{t} \in \Omega.
}
This system is not an orthonormal basis of $\rL^2(\Omega)$, but it is a tight, linearly-independent frame with $A=B=1$.   If we introduce the truncated frames
\be{
\label{FE_indexing}
\Phi_N = \{ \phi_{\bm{n}} \}_{\bm{n} \in I_N},\quad I_N = \left \{ \bm{n} = (n_1,\ldots,n_d) \in \bbZ^d : -\tfrac12 N^{1/d} \leq n_1,\ldots,n_d < \tfrac12 N^{1/d} \right \},
}
(we assume throughout that $\frac12 N^{1/d}$ is a positive integer) then the convergence rate of the resulting orthogonal projections $\cP_N f$ is spectral (Proposition \ref{l:FE_bounded_conv}). 

The approximation based on the frame \R{FE_frame} is known as a Fourier \textit{extension} (or \textit{continuation}) \cite{BoydFourCont,brunoFEP} in the one-dimensional case, and occasionally referred to as a Fourier \textit{embedding} in higher dimensions \cite{boyd2005fourier,pasquettiFourEmbed}.  The connection to frames was first explored in \cite{huybrechs2010fourier}, and further analysis of the $d = 1$ case was given in \cite{BADHFEResolution,FEStability,LyonFESVD}.

Recalling the discussion in \S \ref{ss:canonical_dual_frame}, this frame is an example where the canonical dual frame expansion \R{dual_rep} converges slowly.  Indeed, since $\Phi$ is tight with $A = B = 1$ it is its own canonical dual frame, and therefore the frame coefficients are $a_{\bm{n}} = \ip{f}{\phi_{\bm{n}}}$.  They are precisely the Fourier coefficients of the extension $\tilde{f}$ of $f$ by zero to $(-1,1)^d$:
\bes{
\ip{f}{\phi_{\bm{n}}} =  \int_\Omega f(\bm{t}) \phi_{\bm{n}}(\bm{t}) {\rm d}\bm{t} = \int_{(-1,1)^d} \tilde{f}(\bm{t}) \phi_n(\bm{t}) {\rm d}\bm{t}.
}
As a result, the expansion \R{dual_rep} is nothing more than the Fourier series of $\tilde{f}$ restricted to $\Omega$.  Unless $f$ vanishes smoothly on the boundary $\partial \Omega$, this expansion converges slowly and suffers from Gibbs' phenomenon near $\partial \Omega$. In contrast, the convergence of the best approximation $\cP_N f$ is spectral, regardless of the shape of $\Omega$.

\begin{remark}
This example illustrates a general principle: the restriction of a Riesz basis on a domain $\Omega_e$ to a subset $\Omega \subset \Omega_e$ always results in a frame. If the basis on $\Omega_e$ is orthonormal, the corresponding frame on $\Omega$ is tight. Such a construction (albeit without the connection to frame theory) has been used for the numerical solution of PDEs in complex geometries.  Recent examples of embedding methods implicitly based on such `extension' frames include \cite{lui2009embedding,stein2015immersed}.  See also \cite{albin2011,bruno2010high,lyon2010high} for a method based on one-dimensional extensions.  Besides the Fourier extensions, other examples of extension schemes include volume penalty methods~\cite{shirokoff2015volumepenalty}, immersed domain methods~\cite{boffi2015immersed}, and embedded or fictitious domain methods~\cite{kasolis2015fictitious}.
\end{remark}

\vspace{1pc} \noindent
\textbf{Example 2.\ Augmented Fourier basis.} Consider the case of smooth, nonperiodic functions $f : [-1,1] \rightarrow \bbR$.  Polynomial bases have good convergence properties for such functions but relatively bad resolution power for oscillatory functions.  On the other hand, the Fourier approximation of a nonperiodic function suffers from Gibbs' phenomenon at $t = \pm 1$ and converges only slowly in the $\rL^2$-norm.  One way to seek to remedy this situation is to augment the Fourier basis by a finite number $K \in \bbN$ of polynomials $\psi_1,\ldots,\psi_K$, leading to the system
\be{
\label{AugFourier_frame}
\Phi = \{ \varphi_n \}_{n \in \bbZ} \cup \{ \psi_k \}^{K}_{k=1},\qquad \varphi_n(t)  = \tfrac{1}{\sqrt{2}} \E^{\I n \pi t}.
}
To save unnecessary generalizations, we assume that $\psi_k = \sqrt{k+1/2} P_k$, where $P_k \in \bbP_k$ is the $k^{\rth}$ Legendre polynomial.  Note that $\{ \psi_k \}^{K}_{k=1}$ is an orthonormal basis for
\bes{
\bbP^{0}_K = \left \{ p \in \bbP_K : \int^{1}_{-1} p(t) \D t = 0 \right \}.
}
Since $\{ \varphi_n \}_{n \in \bbZ}$ is an orthonormal basis and $K$ is finite, $\Phi$ forms a frame for $\rH = \rL^2(-1,1)$ with frame bounds $A=1$ and $B=2$. 
It is also linearly independent, since no finite sum of the complex exponentials $\varphi_n$ is exactly equal to a nonconstant algebraic polynomial.  If
\bes{
\Phi_N = \{ \varphi_n : n = -\tfrac{N-K}{2},\ldots,\tfrac{N-K}{2}-1 \} \cup \{ \psi_k \}^{K}_{k=1},\quad N \geq K,\quad \mbox{$N-K$ even},
}
is the truncated frame  (we will not consider the odd case, although it presents few difficulties), then orthogonal projections with respect to this frame inherit the optimal resolution properties of the Fourier basis, yet converge algebraically at rate $K$ for all sufficiently smooth functions (Proposition \ref{l:AugFourier_conv}).  Conversely, the canonical dual frame expansion converges at roughly the same rate as the Fourier expansion of $f$ (see Proposition SM3.2 of the supplementary material \cite{framespart1SM}).

The idea of augmenting the Fourier basis with a finite number of additional functions is an old one, arguably dating back to Krylov \cite{Kryloveappcalc}.  These functions endow the frame with good approximation properties by implicitly subtracting the jump discontinuity of $f$ and its first $K-1$ derivatives at $t = \pm 1$.  This smoothed function now has faster decaying Fourier coefficients, leading to more rapid convergence of $\cP_N f$.  This approach is also referred to as Eckhoff's method \cite{Eckhoff1} or Euler--MacLaurin interpolants \cite{javed2016eulermaclaurin}.  Whilst the convergence of this approximation has been extensively studied (see \cite{BAthesis,BA3} and references therein), the connection with frame theory is, to the best of our knowledge, new.

\begin{remark}
An overall principle illustrated by this example is that adding a finite number of elements of $\rH$ to a Riesz (in particular, orthogonal) basis always results in a frame; a so-called \textit{Riesz frame} \cite[Sec.\ 6.2]{christensen2003introduction}.  Note that such a frame is generally not linearly independent.
Although we focus on polynomials enhancing the Fourier basis here, augmenting a basis with additional terms to incorporate features of the function to be approximated (in this case, smoothness) is quite a general idea.  Other examples include piecewise polynomial functions in the presence of interior discontinuities, or compactly-supported functions in the case of local variations such as oscillations.
\end{remark}

\vspace{1pc} \noindent
\textbf{Example 3.\ Polynomials plus modified polynomials.}
Several problems in numerical analysis call for the approximation of functions of the form
\be{
\label{fw_explicit}
f(t) = w(t) g(t) + h(t),\quad t \in [-1,1],
}
where $g,h$ are smooth functions but $w \in \rL^\infty(-1,1)$ may be singular, oscillatory or possessing some other kind of feature which makes approximation difficult.  The presence of $w(t)$ usually means that the polynomial approximation of $f$ converges only slowly in $N$. One particular instance is
\be{
\label{Shu_frame}
w(t) = (1+t)^{\alpha},\quad 0 < \alpha <1,
}
which corresponds to a weak endpoint singularity of the function $f$.  Note that \R{Shu_frame} has been considered in \cite{ShuEndpoint1,ShuEndpoint2} in the context of PDEs with endpoint singularities.  For other examples corresponding to an oscillatory function $w(t)$, see \cite{hewett2015,GelbPlatteEdge}.

If $w(t)$ is known or has been approximated (as in the applications mentioned above), then it is natural to use it to construct a frame to approximate $f$.  Let $\{ \varphi_n \}_{n \in \bbN}$ be the orthonormal basis of Legendre polynomials (we could also use Chebyshev polynomials here, but we shall use Legendre for simplicity).  Then we form the system
\be{
\label{LegConcat}
\Phi = \{ \varphi_n \}_{n \in \bbN} \cup \{ \psi_n \}_{n \in \bbN},\qquad \psi_n(t) = w(t) \varphi_n(t).
}
Since $w \in \rL^\infty(-1,1)$ and $\{ \varphi_n \}_{n \in \bbN}$ is an orthonormal basis, this system gives rise to a frame for the space $\rH = \rL^2(-1,1)$.  The frame bounds are
\bes{
A = 1 + \essinf_{t \in (-1,1)} | w(t) |^2,\quad B = 1 + \esssup_{t \in (-1,1)} |w(t) |^2.
}
This frame is linearly dependent if and only if $w$ is a rational function of two polynomials.  We assume from now on that this is not the case.  For even $N$, we define the truncated frames by $\Phi_N = \{ \varphi_n \}^{N/2}_{n=1} \cup \{ \psi_n \}^{N/2}_{n=1}$. Orthogonal projections with respect to this frame are spectrally convergent with respect to the smoothness of $g$ and $h$ (Proposition \ref{l:PMP_conv}).  Conversely, the convergence of the dual frame expansion is generally not spectral, but algebraic at a fixed and low rate (Proposition SM4.2 of \cite{framespart1SM}).

Note that more terms can be included in \R{fw_explicit}, i.e.\ $f(t) = \sum^{K}_{i=1} w_i(t) g_i(t)$ for functions $w_1,\ldots,w_K$.  If these are known, then this would lead to the frame construction $\Phi = \cup^{K}_{i=1} \{ \psi_{i,n} \}_{n \in \bbN}$, where $\psi_{i,n}(t) = w_i(t) \varphi_n(t)$.  For simplicity, we consider only the case \R{fw_explicit}, although the generalization is conceptually straightforward. As in Example 2, the interpretation of this approach as a frame approximation has not, to the best to the best of our knowledge, been considered before.

\begin{remark}
The concatenation of a finite number of Riesz or orthonormal bases always results in a frame (albeit one which is generally not linearly independent). More generally, the composition of several frames is still a frame. We note that the concept of concatenations of bases or frames is widely-used in signal and image processing \cite{CandesFramesCS}.  Typically, images and signals may have substantially sparser representations in the resulting frame than in a single orthonormal basis, which yields benefits in tasks such as compression and denoising \cite{mallat09wavelet}.\footnote{For completeness, we note that the system $\{ \psi_n \}_{n \in \bbN}$ in \R{LegConcat} is only a Riesz basis if $w(t)$ is bounded away from zero.  This is not the case in \R{Shu_frame}.  However, $\{ \psi_{n} \}_{n \in \bbN}$ is always a \textit{Bessel sequence} (see, for example, \cite[Def.\ 3.2.2]{christensen2003introduction}); that is, a sequence for which the upper frame condition $\sum_{n \in \bbN} |\ip{f}{\psi_n} |^2 \leq B \| f \|^2$ holds, but for which the lower 
frame condition need not hold.  The composition of a Riesz basis or frame with a finite number of Bessel sequences is also always a frame.}
\end{remark}

\section{Truncated Gram matrices and ill-conditioning}
\label{sect:IllCond}

Since we have assumed linear independence, the truncated frame $\Phi_N$ is also a Riesz basis for its span. Hence the truncated frame bounds $A_N$ and $B_N$ are the same as the Riesz bounds. However, this finite basis is very skewed, and this results in the ill-conditioning of truncated frames. In this section we explore this phenomenon in more detail.

\subsection{The Gram operator}

The Gram operator \R{GramOp} of a frame $\Phi$ is self-adjoint and nonnegative on $\ell^2(I)$ with closed range.  It is bounded, and its restriction $\cG : \ell^2(I) \rightarrow \mathrm{Ran}(\cG)$ is invertible.  Its spectrum $\sigma(\cG)$ satisfies
\bes{
\{0 \} \subseteq \sigma(\cG) \subseteq \{ 0 \} \cup [A,B],
}
where $A,B$ are the frame bounds \cite{HarrisonThesis}.  As shown in \cite{BenedettoTeolisLocalFrame}, $\cG$ is compact if and only if $\rH$ is finite-dimensional, and $\cG$ is positive if and only if $\Phi$ is a Riesz basis.  Hence, in this paper $\cG$ is singular and thus $\mathrm{Ker}(\cG) \neq \{ \bm{0} \}$.  Nonetheless, since $\cG$ has closed range, we may define its Moore--Penrose pseudoinverse $\cG^{\dag} : \ell^2(I) \rightarrow \ell^2(I)$ \cite{BenedettoTeolisLocalFrame}.  One then has the following relation between $\cG$, $\cG^{\dag}$ and the frame bounds:
\be{
\label{GramOpFrameBounds}
A = \nmu{\cG^{\dag}}^{-1},\quad B = \nm{\cG}.
}
Here $\nm{\cdot}$ is the operator norm on $\ell^2(I)$.

\subsection{Truncated Gram matrices}
We now consider the condition number of the matrix $\bm{G}_N$.  Following the discussion above note the following {standard result:}

\lem{
\label{l:GramNormFrameBds}
The truncated Gram matrix $\bm{G}_N$ of a linearly-independent frame $\Phi$ is invertible with $\| \bm{G}^{-1}_N \|^{-1} = A_N$ and $\| \bm{G}_N \| = B_N$,  
where $A_N$ and $B_N$ are the frame bounds of the truncated frame $\Phi_N$. In particular, the condition number of $\bm{G}_N$ is precisely the ratio of the truncated frame bounds: $\kappa(\bm{G}_N)= \| \bm{G}_N \| \| \bm{G}^{-1}_N \| = B_N/ A_N$.
}

In practice, we will also use the following characterization of the frame bounds:
\be{
\label{FrameBds_Char}
A_N = \min_{\substack{\bm{x} \in \bbC^N \\ \| \bm{x} \|=1}} \nm{\cT_N \bm{x} }^2,\quad B_N = \max_{\substack{\bm{x} \in \bbC^N \\ \| \bm{x} \|=1}} \nm{\cT_N \bm{x} }^2,
}
which follows from the fact that $\bm{G}_N = \cT^*_N \cT_N$.
This characterization leads to several different interpretations of the frame bounds. First, $A_N$ measures how small an element of $\rH_N$ can be while having unit-norm coefficients in the frame.  Equivalently, it measures how well the zero element $0 \in \rH$ can be approximated by a nonzero element $\rH_N$ with unit-norm coefficients.
It is clear from this that when frame elements are close to being linearly dependent, the constant $A_N$ can be quite small.
Conversely, $B_N$ measures how large an element $\rH$ with unit-norm coefficients can be. Hence, one expects $B_N$ to remain bounded even for near linearly-dependent frame elements.

\lem{
\label{lem:pollution}
Let $\Phi$ be a linearly-independent frame.  Then
\begin{itemize}
\item[(i)] the sequences $\{ A_N \}_{N \in \bbN}$ and $\{ B_N \}_{N \in \bbN}$ are monotonically nonincreasing and nondecreasing respectively,
\item[(ii)] $B_N \leq B$ for all $N$ and $B_N \rightarrow B$ as $N \rightarrow \infty$,
\item[(iii)] $\inf_{N} A_N >0$ if and only if $\Phi$ is a Riesz basis.
\end{itemize}
}
\prf{
Part (i) follows immediately from \R{indexsets} and \R{FrameBds_Char}, as does the the observation that $B_N \leq B$ in part (ii).  To deduce convergence, let $0 < \epsilon < \sqrt{B}$ be arbitrary and suppose that $\bm{x} \in \ell^2(I)$, $\| \bm{x} \| =1$, is such that $\sqrt{B} \geq \| \cT \bm{x} \| = \sqrt{\ip{\cS \bm{x}}{\bm{x}} }\geq \sqrt{B} - \epsilon$.  Let $\bm{z} \in \bbC^N$ be such that $z_n = x_n$ for $n \in I_N$ and suppose that $\bm{x}_N \in \ell^2(I)$ is the extension of $\bm{z}$ by zero.  Then
\bes{
\sqrt{B_N} \geq \frac{\| \cT_N \bm{z} \|}{\| \bm{z} \|} = \frac{\| \cT \bm{x}_N \|}{\| \bm{x}_N \|} \geq \frac{\| \cT \bm{x} \| - \| \cT ( \bm{x}- \bm{x}_N) \|}{\| \bm{x} \| + \| \bm{x} -\bm{x}_N \|} \geq \frac{\sqrt{B} - \epsilon - \sqrt{B} \| \bm{x} - \bm{x}_N \|}{1 + \| \bm{x} - \bm{x}_N \|}.
}
Since $\bm{x}_N \rightarrow \bm{x}$ as $N \rightarrow \infty$, we deduce that $\sqrt{B} \geq \sqrt{B_N} \geq \sqrt{B} -2 \epsilon$ for all sufficiently large $N$.  Since $\epsilon$ was arbitrary we see that $B_N \rightarrow B$.  Finally, for part (iii) we use \cite[Prop.\ 6.1.2]{christensen2003introduction}.
}

These lemmas imply that the truncated Gram matrices $\bm{G}_N$ are necessarily ill-conditioned for large $N$.  Such ill-conditioning can also be arbitrarily bad:

\prop{
\label{p:simple_bad_frame}
Let $\{ \varphi_n \}_{n \in \bbN}$ be an orthonormal basis of $\rH$ and let $g \in \rH$, $\| g \| =1$, be such that $\ip{g}{\varphi_n} \neq 0$ for infinitely many $n$.  Then the system $\Phi = \{ g , \varphi_1,\varphi_2,\ldots \}$ is a linearly-independent frame for $\rH$ with bounds $A = 1$ and $B=2$.  Moreover, if $\Phi_N = \{ g , \varphi_1,\phi_2,\ldots,\varphi_{N-1} \}$, then the finite frame bounds are given by
\be{
\label{simple_exact}
A_N = 1 - \sqrt{\sum^{N-1}_{n=1} | \ip{g}{\varphi_n} |^2},\qquad B_N = 1 + \sqrt{\sum^{N-1}_{n=1} | \ip{g}{\varphi_n} |^2}.
}
}

See \S SM1 of the supplementary material for the proof. Observe that $A_N \leq \sqrt{\sum_{n \geq N} | \ip{g}{\varphi_n} |^2 }$ since $\| g \| =1$ for this frame.  Hence the behaviour of $A_N$ is related to the decay of the coefficients $\ip{g}{\varphi_n}$ in the basis $\{ \varphi_n \}_{n \in \bbN}$.  The better $g$ is approximated in this basis, the worse the conditioning of the truncated Gram matrix.  

This result illustrates how easy it is for the truncated Gram matrices to be ill-conditioned: we can create arbitrarily bad conditioning merely by augmenting an orthonormal basis with one additional element that is well approximated in the basis.  When more elements are added (as in Example 2) or a whole orthonormal basis is added (as in Example 3), it is therefore not surprising that the corresponding Gram matrices can be exceedingly ill-conditioned.  See \S \ref{ss:example_condition_numbers}.

\begin{figure}
\begin{center}
$\begin{array}{ccc}
\includegraphics[width=5.00cm]{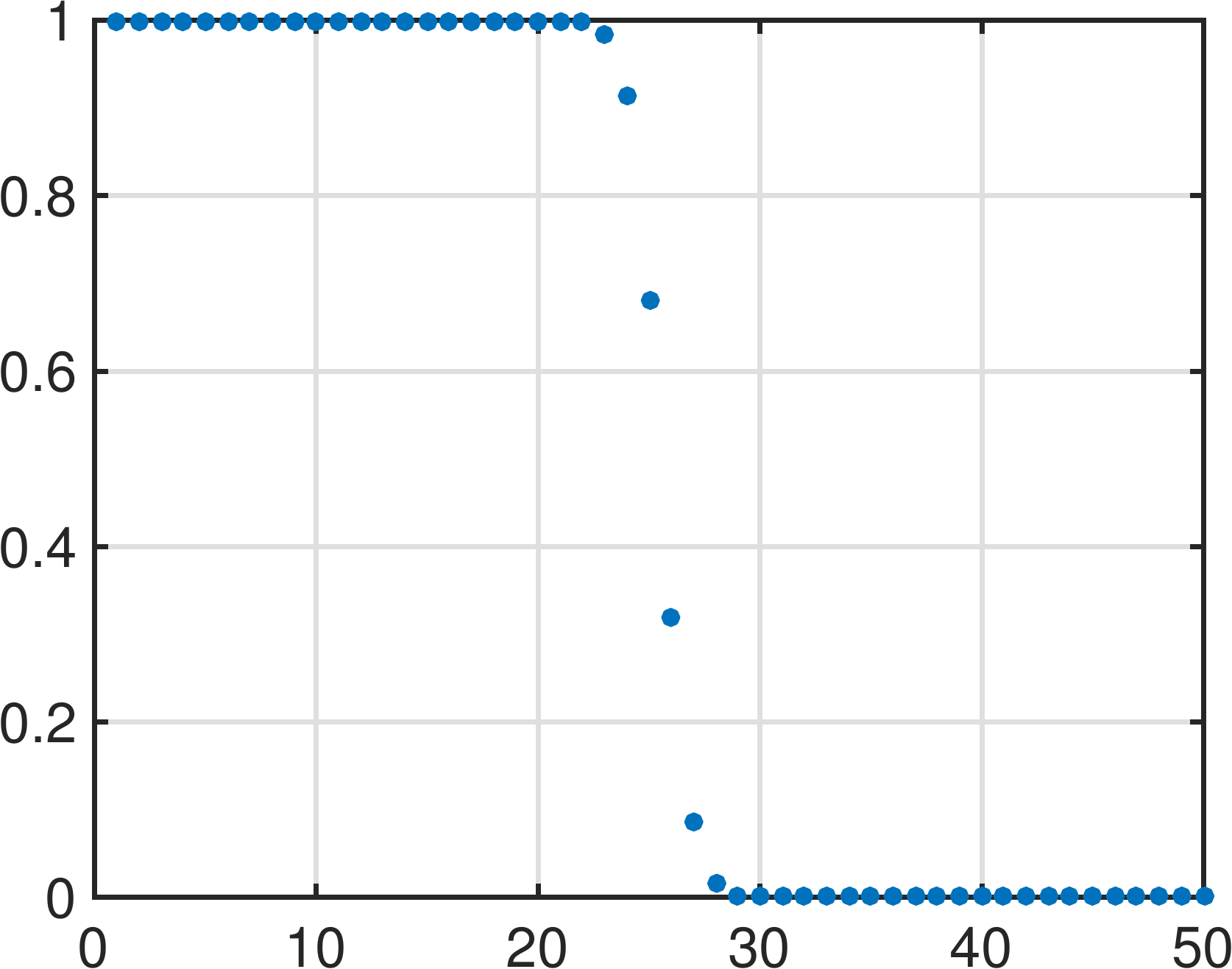} 
&\includegraphics[width=5.00cm]{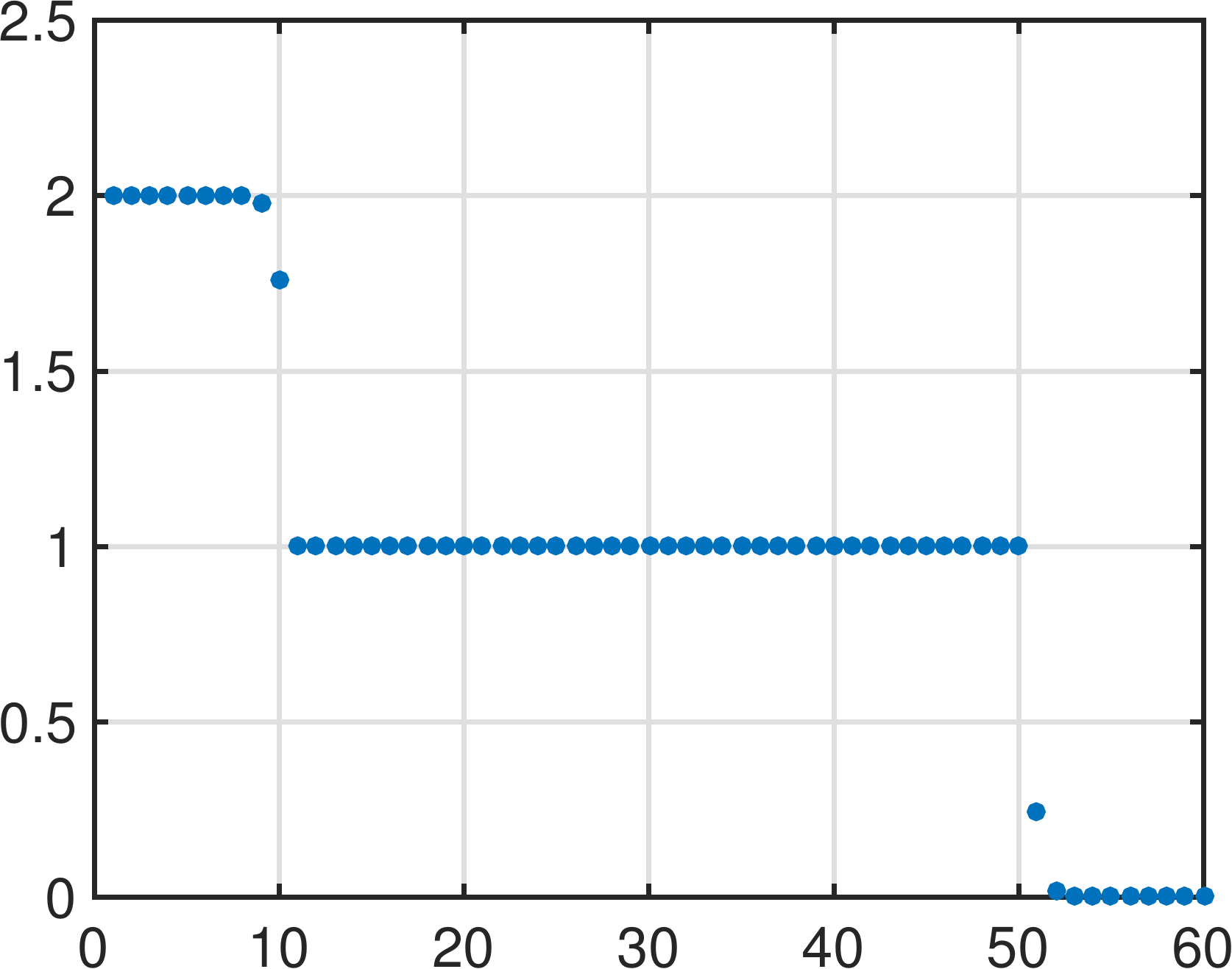}  & \includegraphics[width=5.00cm]{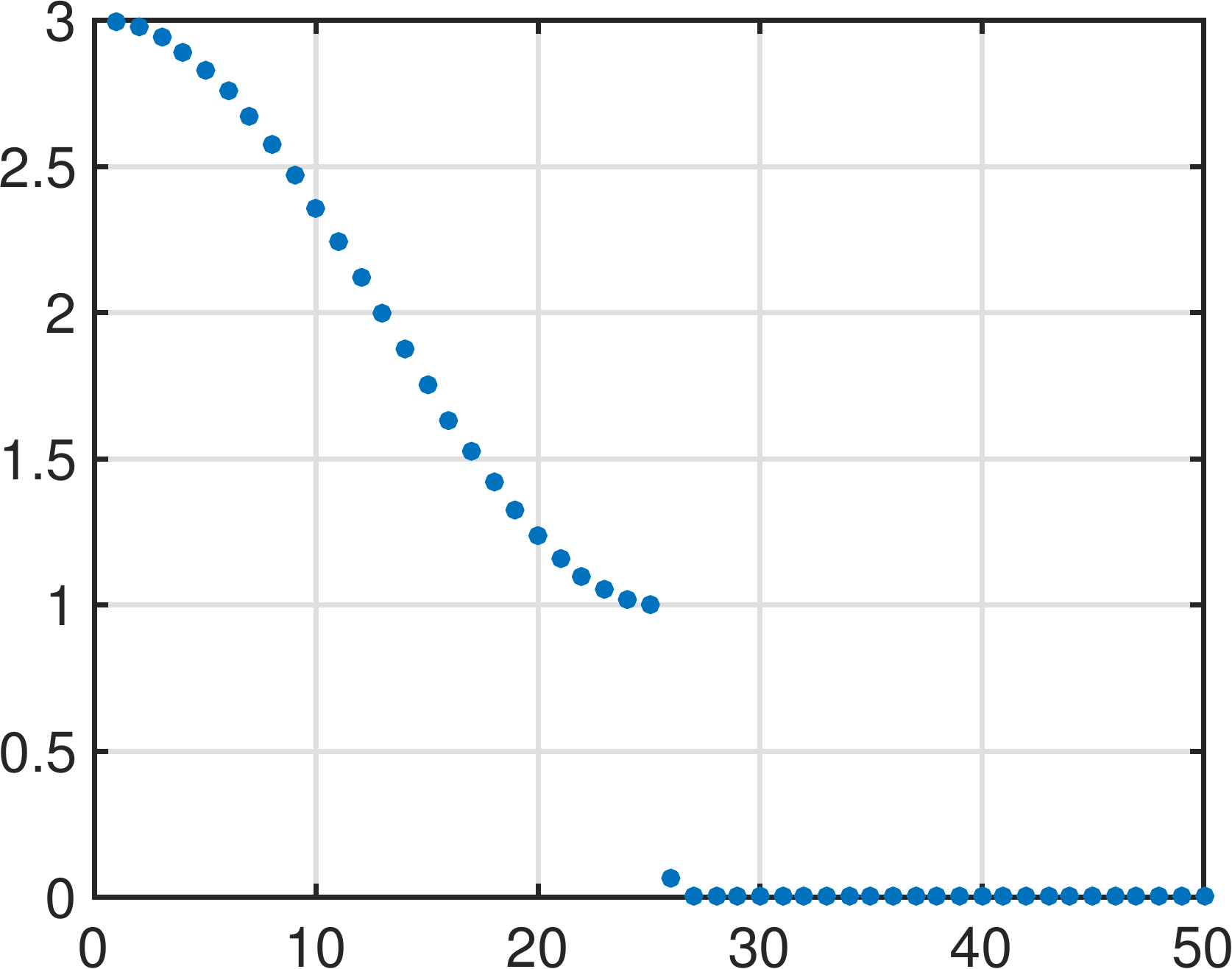} 
\end{array}$
\caption{The eigenvalues of $\bm{G}_{50}$ for Examples 1--3 (left to right).  The parameters used were $T=2$ (Example 1), $K=10$ (Example 2) and $w(t) = \sqrt{1+t}$ (Example 3).} 
\label{f:GramSpectrum}
\end{center}
\end{figure}

\rem{
\label{r:opthy}
The truncated Gram matrices $\bm{G}_N$ are the $N \times N$ finite sections of the Gram operator $\cG$.  Recall that $\{ 0\} \subseteq \sigma(\cG) \subseteq \{0\} \cup [A,B]$.  Unfortunately, as illustrated in Fig.\ \ref{f:GramSpectrum} for the three examples, $\sigma(\bm{G}_N)$ does not lie within $\{ 0 \} \cup [A,B]$.  Instead, spurious eigenvalues are introduced in the spectral gap between $0$ and the lower frame bound $A$; a well-known phenomenon referred to as pollution in the finite section method \cite{DaviesPlumSpecPoll,LewinSereSpecPoll}.  The small eigenvalues of $\bm{G}_N$ are of particular relevance to this paper.  These translate into ill-conditioning of $\bm{G}_N$, meaning that the system \R{GramLinSyst} requires regularization.  Note that $0 \notin \sigma(\bm{G}_N)$ since the frame is linearly independent, yet as seen in Lemma \ref{lem:pollution}, small eigenvalues necessarily arise as approximations to the zero eigenvalue of $\cG$. 
}

\subsection{Examples}
\label{ss:example_condition_numbers}

We now consider $\kappa(\bm{G}_N)$ for the examples of \S \ref{sect:examples}. Proofs of the results presented in this section can be found in the supplementary material \cite{framespart1SM}.

\begin{figure}
\begin{center}
$\begin{array}{ccc}
\includegraphics[width=5.00cm]{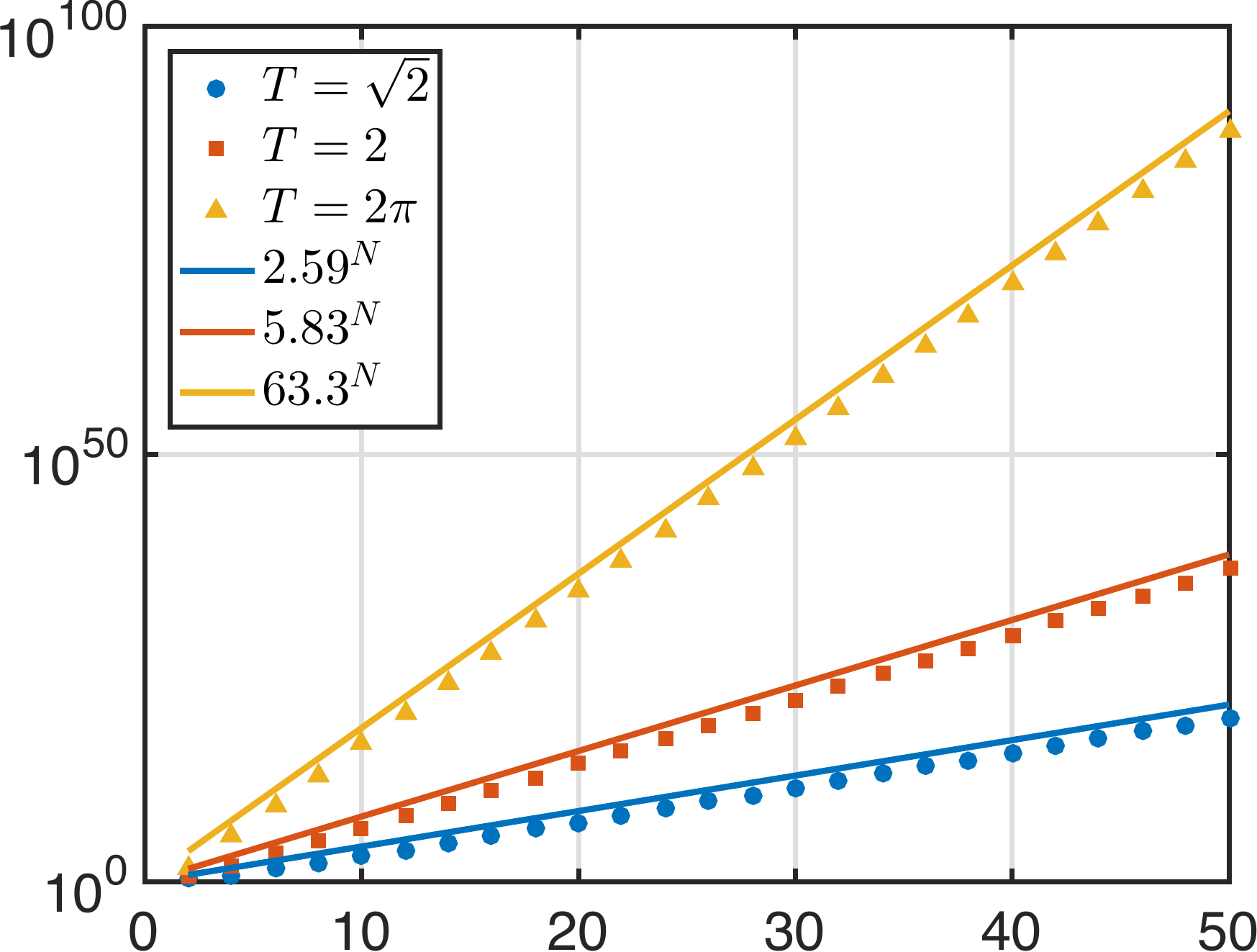}  
&  \includegraphics[width=5.00cm]{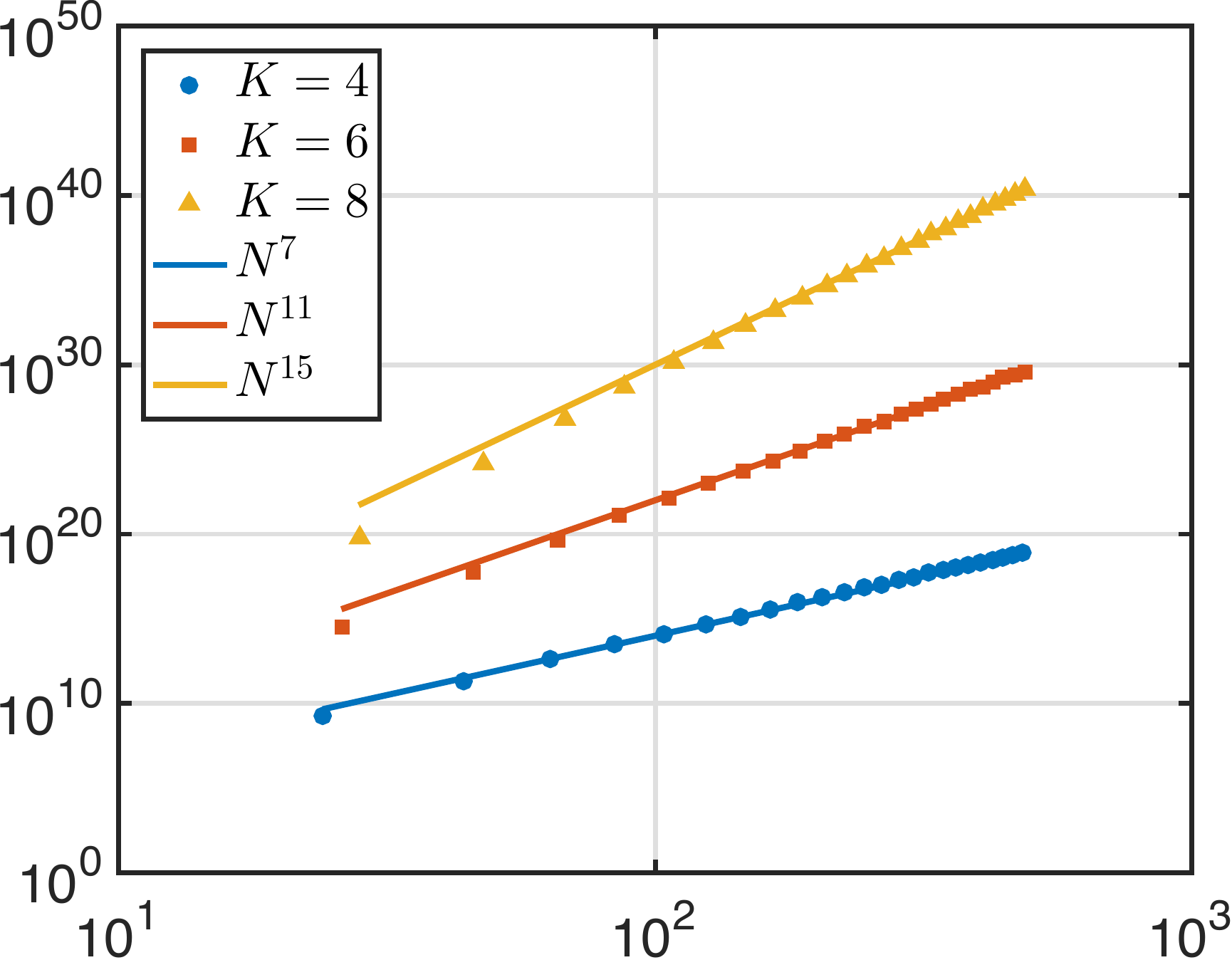}  &
\includegraphics[width=5.00cm]{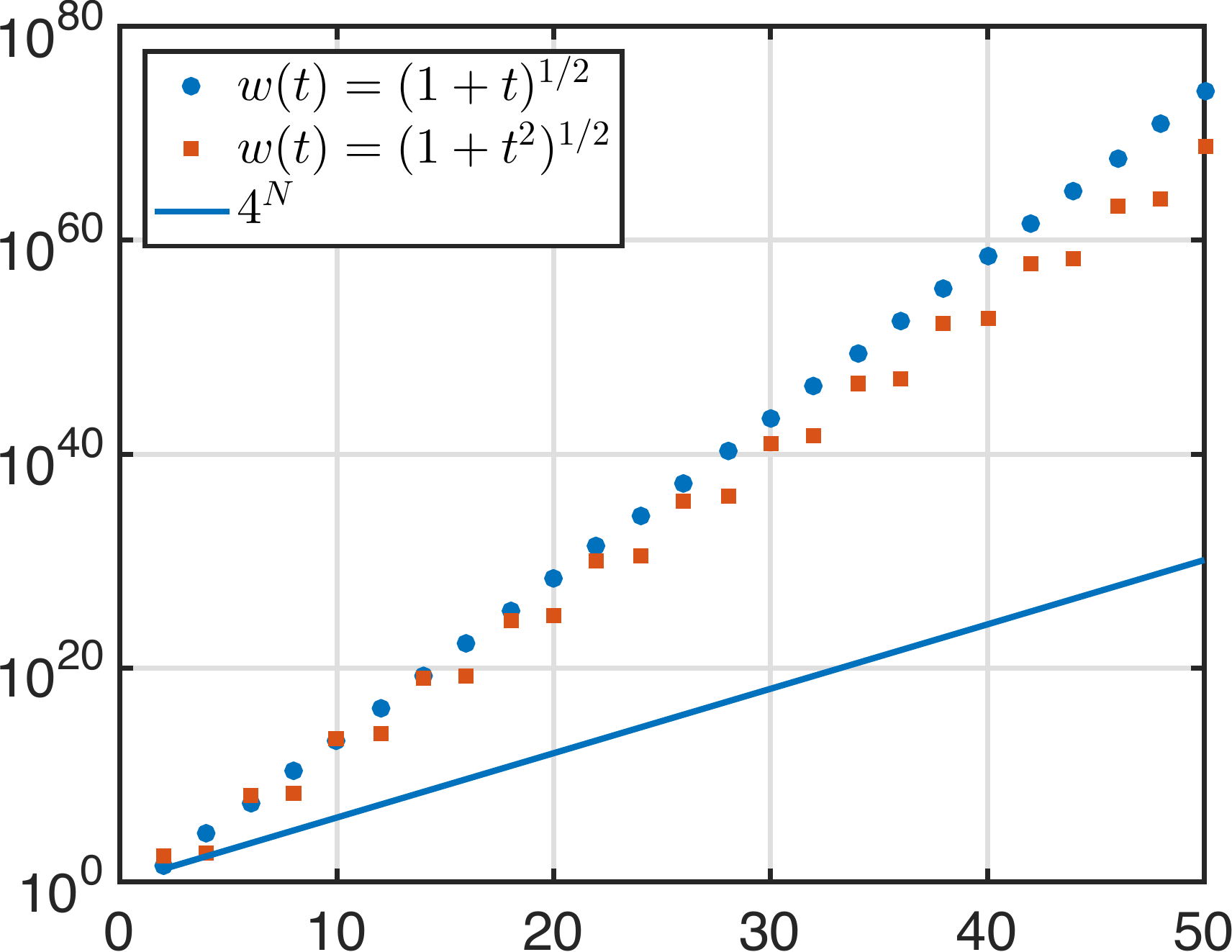} 
\end{array}$
\caption{The condition numbers $\kappa(\bm{G}_N)$ for Examples 1--3 (left to right) with parameters $T = \sqrt{2},2,2\pi$, $K=4,6,8$ and $w(t) = \sqrt{1+t}$  respectively.  The solid lines show the bounds in Propositions \ref{l:FE_illcond}--\ref{l:LegConcat_Gram}.  Computations were carried out using additional precision.} 
\label{f:GramCondition}
\end{center}
\end{figure}

\vspace{1pc}\noindent
\textbf{Fourier frames for complex geometries.}  Consider the frame of Example 1.  If $\Omega = (-\frac1T,\frac1T)$ is an interval, where $T>1$, then $\kappa(\bm{G}_N) = \ord{E(T)^{N}}$ as $N \rightarrow \infty$, where $E(T) = \cot^2(\pi/(4T)) > 1$ \cite{FEStability}.  Hence the condition number is geometrically large in $N$ -- see Fig.\ \ref{f:GramCondition}(a).  A similar, albeit somewhat weaker, result also holds in arbitrary dimensions:

\prop{
\label{l:FE_illcond}
Let $\Omega \subseteq (-1,1)^d$ be a compact, Lipschitz domain and consider the frame \R{FE_frame}.  Then the condition numbers $\kappa(\bm{G}_N)$ grow superalgebraically fast in $N$.
}

The explanation of this result is rather simple.  One can show that the kernel $\mathrm{Ker}(\cG)$, a subset of $\ell^2(I)$, consists precisely of the sequences of Fourier coefficients of those functions defined on $(-1,1)^d$ which vanish on $\Omega$ (Proposition SM2.1 of the supplementary material). 
Now consider a smooth function $g$ with this property.  Then its Fourier expansion on $(-1,1)^d$ converges superalgebraically fast to $g$ on $(-1,1)^d$, and therefore to zero on the domain $\Omega$.  The ratio of the norm of the truncated Fourier series on $\Omega$ divided by its norm on $(-1,1)^d$ is an upper bound for $A_N$, implying ill-conditioning at a superalgebraic rate.

\vspace{1pc}\noindent
\textbf{Augmented Fourier basis.} Consider the frame of Example 2. 
In this case, $\kappa(\bm{G}_N)$ grows algebraically fast at a rate depending on $K$:

\prop{
\label{l:AugFourier_Gram}
Let $K \in \bbN$ be fixed and consider the frame \R{AugFourier_frame}.  Then $\kappa(\bm{G}_N) \gtrsim N^{2K-1}$ as $N \rightarrow \infty$.
}

The intuition behind this result is as follows.  The kernel of the Gram operator $\cG$ has dimension $K$, and consists of sequences which are comprised of the coefficients $\{ \ip{p}{\psi_k} \}^{K}_{k=1}$ and $\{ \ip{-p}{\varphi_n} \}_{n \in \bbZ}$, where $p$ is an arbitrary polynomial in $\bbP^0_{K}$.  See Proposition SM3.1.  It is possible to construct a polynomial $p \in \bbP^0_K$ which has $K$ orders of periodic smoothness (of course, $p$ is analytic, but it is not periodic in general).  This function has Fourier coefficients $\ip{p}{\varphi_n} $ which decay like $|n|^{-k-1}$ as $n \rightarrow \pm \infty$.  Hence there is a function in $\rH_N$, i.e.\ the difference between $p$ and its partial Fourier series, which is of magnitude $\ord{N^{-K}}$ but which has $\ord{1}$ coefficients in the frame $\Phi_N$.

\vspace{1pc}\noindent
\textbf{Polynomials plus modified polynomials.} For the frame of Example 3 we have:

\prop{
\label{l:LegConcat_Gram}
Let $\Phi$ be the frame \R{LegConcat} with $w(t)$ given by \R{Shu_frame}.  Then $\kappa(\bm{G}_N) \gtrsim 4^{N}$ as $N \rightarrow \infty$ up to an algebraic factor in $N$.  
}

The idea behind this result is similar to that of Example 2 (see Proposition SM4.1 for a description of $\mathrm{Ker}(\cG)$ in this case). We choose a polynomial $q(t) = (1+t)^{N/2-1}$ such that $w(t)q(t)$ has several orders of smoothness at $t=0$ in spite of the algebraic singularity there. Hence, it can be well approximated by a single polynomial $p(t)$. The difference $p-wq \in \rH_N$ is close to zero but it has $\ord{1}$ coefficients in the frame.

\vspace{1pc}\noindent
Numerical illustrations of these three estimates are shown in Fig.\ \ref{f:GramCondition}.  Unlike in Examples 1 and 2, the lower bound of $4^N$ in Example 3 does not give a good estimate of the true growth of $\kappa(\bm{G}_N)$.

\section{Computing best approximations}
\label{sect:OP_compute}

We now turn our attention to the computation of the orthogonal projection $\cP_N f$. Our first assertion, developed in the next subsection, is that it is impossible in general to compute $\cP_N f$, since the coefficients of this approximation can grow arbitrarily fast with $N$. 

%However, with a simple numerical scheme one can compute best approximations subject to having small norm coefficients. This numerical solution is computable, and more desirable in practice since it results in bounded coefficients.

\subsection{Impossibility of computing best approximations}\label{ss:impossible}

Computing $\cP_N f = \sum_{n \in I_N} x_n \phi_n$ requires solving the ill-conditioned linear system \R{GramLinSyst}.  If $\bm{x} = \{ x_n \}_{n \in I_N}$ and $\bm{y} = \{ \ip{f}{\phi_n} \}_{n \in I_N}$ then by Lemma \ref{l:GramNormFrameBds}, the Cauchy--Schwarz inequality and the frame condition \R{framecond} we have
\bes{
\| \bm{x} \| = \| \bm{G}^{-1}_{N} \bm{y} \| \leq \| \bm{G}^{-1}_N \| \| \bm{y} \| \leq A^{-1}_N \sqrt{\sum_{n \in I} | \ip{f}{\phi_n} |^2 }  \leq A^{-1}_{N} \sqrt{B} \| f \|.
}
Hence, the coefficients $\bm{x}$ of the orthogonal projection may, in the worst case, grow as rapidly as $A^{-1}_N$.  Of course, this is only an upper bound, and therefore may not be achieved for a fixed $f \in \rH$.  However, it is easy to create an example where the growth of $\|\bm{x} \|$ mirrors that of $A^{-1}_N$. 

\prop{
\label{p:unbounded_solution}
Let $\{ \varphi_n \}_{n \in \bbN}$, $g \in \rH$, $\Phi$ and $\Phi_N$ be as in Proposition~\ref{p:simple_bad_frame} and suppose that $f \in \rH$, $\| f \|=1$ is given by
\bes{
f = \frac{\sqrt{6}}{\pi} \sum_{n \in \bbN} \frac{\mathrm{sign}(\ip{g}{\varphi_n})}{n} \phi_n,
}
where, for $\omega \in \bbC$, $\mathrm{sign}(\omega) = \omega / |\omega|$ if $\omega \neq 0$ and $\mathrm{sign}(\omega)=0$ otherwise.  Suppose also that $\sup_{n \in \bbN} |\ip{g}{\varphi_n}| n^2 < \infty$.  If $\bm{x} = \{ x_n \}_{n \in I_N}$ is the solution of \R{GramLinSyst} then
\bes{
\| \bm{x} \| \geq \left ( \frac{\pi}{\sqrt{6}} \max_{n \geq N} \left \{ n | \ip{g}{\varphi_n} | \right \} \right )^{-1}.
}
}
%\prf{
%The orthogonality of the basis $\{ \phi_n \}_{n \in \bbN}$ and the fact that $\| g \| =1$ means that the system $\bm{G}_N x = y$ is equivalent to
%\eas{
%x_0 + \sum^{N-1}_{m=1} \overline{\ip{g}{\phi_m}} x_{m} = \ip{f}{g},\qquad \ip{g}{\phi_n} x_0 + x_{n}  = \ip{f}{\phi_n},\quad n=1,\ldots,N-1.
%}
%Substituting the second equation into the first gives
%\bes{
%x_0 \left ( 1 - \sum^{N-1}_{m=1} | \ip{g}{\phi_m} |^2 \right ) + \sum^{N-1}_{m=1} \ip{f}{\phi_m}\overline{\ip{g}{\phi_m}} = \ip{f}{g},
%}
%and after rearranging and using the definition of $f$, this gives
%\be{
%\label{x0express}
%x_0 = \frac{\sum_{m \geq N} \ip{f}{\phi_m}\overline{\ip{g}{\phi_m}} }{\sum_{m \geq N} | \ip{g}{\phi_m} |^2} = \frac{\sqrt{6} \sum_{m \geq N} n^{-1} | \ip{g}{\phi_m} |}{\pi \sum_{m \geq N} | \ip{g}{\phi_m} |^2}.
%}
%Observe that
%\bes{
%\sum_{m \geq N} | \ip{g}{\phi_m} |^2 \leq \max_{n \geq N} \left \{ n | \ip{g}{\phi_n} | \right \} \sum_{m \geq N} m^{-1} | \ip{g}{\phi_m} |.
%}
%Substituting this into \R{x0express} and noting that $\| x \| \geq |x_0|$ now gives the result.
%}

{See \S SM1 for a proof.}
Suppose now we augment a basis $\{ \varphi_n \}$ with an element $g \in \rH$ that is well approximated in the basis; for example, $|\ip{g}{\varphi_n}| = \ord{n^{-\alpha-1}}$ for some $\alpha \geq 1$.  Then $A^{-1}_{N} \gtrsim N^{\alpha+1/2}$ by Proposition \ref{p:simple_bad_frame} and, if $\bm{x}$ is as in Proposition \ref{p:unbounded_solution}, then $\| \bm{x} \| \gtrsim N^{\alpha}$.  Hence, there exists a fixed $f \in \rH$ whose projection coefficients coefficients $\bm{x}$ grow rapidly as the Gram matrix $\bm{G}_N$ becomes more ill-conditioned.

Although this example is synthetic, it illustrates the general principle that the coefficients of the orthogonal projection $\cP_N f$ in a truncated frame approximation can grow at a similar rate to that of the condition number.  Hence it is generally impossible to compute $\cP_N f$ accurately in floating-point arithmetic.  To see this in a more practical setting, in Table \ref{t:FEcoeffics} we display the coefficients for Example 1 when applied to several different functions.  As is evident, only for the entire function $f(x) = \exp(x)$ is the growth of the coefficients avoided.  For the other two functions, which are less smooth, we witness geometric growth of the coefficients, mirroring that of the condition number (see Proposition \ref{l:FE_illcond}).

\begin{table}[h]
\begin{center}
\begin{tabular}{|c|ccccc|}
            \hline
       $N$     & $10$ & $20$ & $40$ & $80$ & $160$ \\
            \hline
            $f(t) = \exp(t)$ &   1.77e0 & 1.81e0 & 1.84e0 & 1.86e0 &  1.87e0 \\
            $f(t) = \frac{1}{1+16 t^2}$ & 2.27e0 & 5.05e1 & 3.64e4 & 2.32e10 & 1.13e22  \\
            $f(t) = |t|^5$ & 2.12e-1 & 3.67e-1 & 1.76e4 & 7.62e26 &  6.09e91  \\
            $\kappa(\bm{G}_N)$ & 1.84e6 & 5.64e13 & 8.01e28 & 2.35e59 & 2.90e120 \\
            \hline
        \end{tabular}
\caption{The $\ell^2$-norm $\| \bm{x} \|$ of the coefficients of the orthogonal projection $\cP_N f$ for Example 1.}\label{t:FEcoeffics}
\end{center}
\end{table}

\subsection{Truncated SVD projections}
With this in mind, we now turn our attention to computing approximate projections. 

To regularize \R{GramLinSyst}, we resort to a familiar approach for treating ill-conditioned and ill-posed problems (see, for example, \cite{EnglRegularization,HansenEtAlLeastSquares,NeumaierIllCond}): compute the SVD of $\bm{G}_N$, discard all singular values below a tolerance $\epsilon$ and then find the solution $\bm{x}^{\epsilon}$ of the resulting system\footnote{We note in passing that Tikhonov regularization or similar could be used instead, with a number of changes to the results that follow.}.  The entries of $\bm{x}^{\epsilon}$ are no longer the coefficients of $\cP_N f$ but rather the projection onto a smaller space $\rH^{\epsilon}_N$ depending on $\epsilon$.  However, we shall show the the following.  First, $\nmu{\bm{x}^{\epsilon}}$ is bounded (by virtue of the truncated SVD), unlike the coefficients of the orthogonal projection $\cP_N f$.  Second, despite many of the singular values having been discarded, the corresponding projection can still approximate $f$ to high accuracy.

We first require some notation.  Since $\bm{G}_N$ is positive definite its singular values $\sigma_1,\ldots,\sigma_N$ are its eigenvalues and its SVD takes the form
\[
\bm{G}_N = \bm{V} \bm{\Sigma} \bm{V}^*,
\]
where $\bm{V} \in \bbC^{N\times N}$ is unitary and $\bm{\Sigma} = \mathrm{diag}(\sigma_1,\ldots,\sigma_N)$ is diagonal.  Write $\{ \bm{v}_n : n \in I_N \}$ for the columns of $\bm{V}$, which are the left/right singular vectors of $\bm{G}_N$, i.e.\ $\bm{G}_N \bm{v}_n = \sigma_n \bm{v}_n$, $n \in I_N$.  To each singular vector we associate a element $\xi_n$ of $\rH_N$:
\bes{
\xi_n = \sum_{m \in I_N} (\bm{v}_n)_m \phi_m = \cT_N \bm{v}_n \in \rH_N.
}
It follows from the orthogonality of $\bm{v}_n$ that the functions $\xi_n$ are orthogonal in $\rH$:
\be{
\label{xi_orthog}
\ip{\xi_n}{\xi_m} = \ip{\cT_N \bm{v}_n}{\cT_N \bm{v}_m} = \ip{\bm{v}_n}{\cT^*_N \cT_N \bm{v}_m} = \sigma_m \ip{\bm{v}_n}{\bm{v}_m} = \sigma_m \delta_{n,m}. 
}
As a result they form an orthogonal basis for $\rH_N$.

Given a tolerance $\epsilon > 0$, let $\bm{\Sigma}^{\epsilon}$ be the diagonal matrix with $n^{\rth}$ entry $\sigma_n$ if $\sigma_n > \epsilon$ and zero otherwise, and define
\be{
\label{trunc_Gram_eps}
\bm{G}^{\epsilon}_N = \bm{V} \bm{\Sigma}^{\epsilon} \bm{V}^*.
}
Then the truncated SVD coefficients $\bm{x}^{\epsilon}$ are $\bm{x}^{\epsilon} = (\bm{G}^{\epsilon}_N)^{\dag} \bm{y} = \bm{V} (\bm{\Sigma}^{\epsilon})^{\dag} \bm{V}^* \bm{y}$,
where $\dag$ denotes the pseudoinverse. Note that $(\bm{\Sigma}^{\epsilon})^{\dag}$ is diagonal with $n^{\rth}$ entry equal to $1/\sigma_n$ if $\sigma_n > \epsilon$ and $0$ otherwise.  We may also write $\bm{x}$ and $\bm{x}^{\epsilon}$ as follows:
\be{
\label{x_singvec}
\bm{x} = \sum_{n \in I_N} \frac{\ip{\bm{y}}{\bm{v}_n}}{\sigma_n} \bm{v}_n,\qquad \bm{x}^{\epsilon} = \sum_{\sigma_n > \epsilon} \frac{\ip{\bm{y}}{\bm{v}_n}}{\sigma_n} \bm{v}_n,
}
where $\ip{\cdot}{\cdot}$ denotes the Euclidean inner product on $\bbC^N$.  Given $\bm{x}^{\epsilon}$, much like with the coefficients $\bm{x}$, we define the approximation
\bes{
\cP^{\epsilon}_{n} f = \cT_N \bm{x}^{\epsilon} = \sum_{n \in I_N} (\bm{x}^\epsilon)_n \phi_n.
}
Observe that $\ip{\bm{y}}{\bm{v}_n} = \sum_{m \in I_N} \ip{f}{\phi_m} \overline{(\bm{v}_n)_m} =\ip{f}{\xi_n}$, and therefore
\be{
\label{Peps_def}
\cP_N f = \sum_{n \in I_N} \frac{\ip{f}{\xi_n}}{\sigma_n} \xi_n, \qquad \cP^{\epsilon}_N f = \sum_{\sigma_{n} > \epsilon} \frac{\ip{f}{\xi_n}}{\sigma_n} \xi_n.
} 
Thus, $\cP^{\epsilon}_N$ is the orthogonal projection from $\rH$ to $\rH^{\epsilon}_N$, where $\rH^{\epsilon}_N = \spn \{ \xi_n : \sigma_n > \epsilon \}$.

\rem{
\label{r:meaning_of_xi}
\label{r:prolates}
Note  that $\| \xi_n \| = \sqrt{\sigma_n}$ due to \R{xi_orthog}. Hence a singular vector $\bm{v}_n$ with small singular value $\sigma_n$ corresponds to a function $\xi_n$ with small norm in $\rH$. In other words, these functions do not contribute much to the approximation. It is therefore not all that surprising that they can be discarded without overly sacrificing accuracy.

In Example 1, the singular vectors $\bm{v}_n$ and the functions $\xi_n$ correspond precisely to the so-called \emph{prolate spheroidal wave sequences} and \emph{prolate spheroidal wave functions}, introduced by Slepian, Landau and Pollak in the study of bandlimited extrapolation \cite{ProlatesIII,SlepianV}. These are a central object of study in the subfield of harmonic analysis and signal processing that focuses on time-frequency localization of signals \cite{daubechies1992tenlectures,HoganLakeyBook}.  
In our setting, the prolate functions corresponding to small singular values are small on $\Omega$ but large on the extended region $(-1,1)^d \backslash \Omega$, i.e.\ they are approximately supported away from $\Omega$ and hence they do not influence the approximation substantially \cite{matthysen2015fastfe}.
}

\subsection{Analysis of the truncated SVD projection $\cP^{\epsilon}_N$}

We now consider the error of the projection $\cP^{\epsilon}_N$.  The following is our main result:
\thm{
\label{t:Proj_err}
The truncated SVD projection $\cP^{\epsilon}_N$ satisfies
\be{
\label{Proj_err}
\| f - \cP^{\epsilon}_N f \| \leq \inf \left \{ \| f - \cT_N \bm{z} \| + \sqrt{\epsilon} \| \bm{z} \| : \bm{z} \in \bbC^N \right \}.
}
}
\prf{
Let $\bm{z} \in \bbC^N$.  Since $\cP^{\epsilon}_N$ is the orthogonal projection onto $\rH^{\epsilon}_N$, we have
\bes{
\| f - \cP^{\epsilon}_N f \| \leq \| f - \cP^{\epsilon}_N \cT_N \bm{z} \| \leq \| f - \cT_N \bm{z} \| + \| \cT_N \bm{z} - \cP^{\epsilon}_N \cT_N \bm{z} \|.
}
Note that $\cT_N \bm{z} = \cP_N \cT_N \bm{z}$ since $\cT_N \bm{z} \in \rH_N$.  Hence \R{Peps_def} and the orthogonality of the $\xi_n$'s gives
\bes{
\| \cT_N \bm{z} - \cP^{\epsilon}_N \cT_N \bm{z} \|^2 = \nm{\sum_{\sigma_n \leq \epsilon} \frac{\ip{\cT_N \bm{z}}{\xi_n}}{\sigma_n} \xi_n}^2 = \sum_{\sigma_n \leq \epsilon} \frac{|\ip{\cT_N \bm{z}}{\xi_n}|^2}{\sigma_n}.
}
Observe that  $\ip{\cT_N \bm{z}}{\xi_n} = \ip{\cT_N \bm{z}}{\cT_N \bm{v}_n} = \ip{\bm{z}}{\bm{G}_N \bm{v}_n} = \sigma_n \ip{\bm{z}}{\bm{v}_n}$ and therefore
\bes{
\| \cT_N \bm{z} - \cP^{\epsilon}_N \cT_N \bm{z} \|^2 =  \sum_{\sigma_n \leq \epsilon} \sigma_n |\ip{\bm{z}}{\bm{v}_n}|^2 \leq \epsilon \sum_{n \in I_N} | \ip{\bm{z}}{\bm{v}_n} |^2 = \epsilon \| \bm{z} \|^2,
}
where in the last step we use the fact that the vectors $\{ \bm{v}_n \}_{n \in I_N}$ are orthonormal.
}

This theorem establishes the claim made earlier in the paper: the convergence of the projection $\cP^{\epsilon}_N f$ is dictated by how well $f$ can be approximated by coefficients $\bm{z}$ with small norm.  Note that this situation is markedly different to the case of the projection $\cP_N f$, wherein the analogous error bound is simply $\| f - \cP_N f \| \leq \inf \{ \| f - \cT_N \bm{z} \| :  \bm{z} \in \bbC^N\}$.  As we discuss in \S \ref{ss:example_convergence}, the appearance of the term $\sqrt{\epsilon} \| \bm{z} \|$ can change the behaviour of the error in a key way.  
We note in passing however that if $\bm{a} =  \{ \ip{f}{\cS^{-1}\phi_n} \}_{n \in I}$ are the frame coefficients of $f$, then setting $\bm{z} = \{ a_n \}_{n \in I_N}$ in \R{Proj_err} and allowing $N \rightarrow \infty$ gives
\be{
\limsup_{N \rightarrow \infty} \nmu{f - \cP^{\epsilon}_N f} \leq \sqrt{\epsilon} \nm{\bm{a}} \leq \sqrt{\epsilon/A} \nm{f},
}
where the second inequality is due to \R{dual_bd}.  Hence the projection $\cP^{\epsilon}_N f$ eventually approximates $f$ to within a factor of $\sqrt{\epsilon}$.

We next consider the behaviour of the coefficients $\bm{x}^{\epsilon}$:

\thm{
\label{t:coefficients1}
The coefficients $\bm{x}^{\epsilon}$ of the truncated SVD projection $\cP^{\epsilon}_N$ satisfy
\be{
\label{coeffics_bounded}
\| \bm{x}^{\epsilon} \| \leq \inf \left \{ 1/\sqrt{\epsilon} \| f - \cT_N \bm{z} \| + \| \bm{z} \| : \bm{z} \in \bbC^N \right \}.
}
Moreover, if $\bm{a}^{\epsilon}_{N} \in \ell^2(I)$ is the extension of $\bm{x}^\epsilon$ by zero, then
\be{
\| \bm{a} - \bm{a}^{\epsilon}_N \| \leq \left ( 1 + \sqrt{B/\epsilon} \right ) \sqrt{\sum_{n \in I \backslash I_N} | a_n |^2} + \sqrt{\epsilon/A} \| \bm{a} \|,
}
where $\bm{a} =  \{ \ip{f}{\cS^{-1}\phi_n} \}_{n \in I}$ are the frame coefficients of $f$. 
}
\prf{
For the first part, we use \R{x_singvec} to write
\bes{
\bm{x}^{\epsilon} = \sum_{\sigma_n > \epsilon} \frac{\ip{f}{\xi_n}}{\sigma_n} \bm{v}_n =\sum_{\sigma_n > \epsilon} \frac{\ip{f-\cT_N \bm{z}}{\xi_n}}{\sigma_n} \bm{v}_n + \sum_{\sigma_n > \epsilon} \frac{\ip{\cT_N \bm{z}}{\xi_n}}{\sigma_n} \bm{v}_n.
}
Consider the first term on the right-hand side.  By \R{xi_orthog} and \R{Peps_def} we have
\bes{
\nm{\sum_{\sigma_n > \epsilon} \frac{\ip{f-\cT_N \bm{z}}{\xi_n}}{\sigma_n} \bm{v}_n}^2 = \sum_{\sigma_n > \epsilon}  \frac{|\ip{f-\cT_N \bm{z}}{\xi_n}|^2}{\sigma^2_n} \leq \frac{1}{\epsilon} \| \cP^{\epsilon}_N (f-\cT_N \bm{z}) \|^2 \leq \frac{1}{\epsilon} \| f - \cT_N \bm{z} \|^2.
}
For the second term, we notice that $\ip{\cT_N \bm{z}}{\xi_n} = \sigma_n \ip{\bm{z}}{\bm{v}_n}$, and therefore
\bes{
\nm{\sum_{\sigma_n > \epsilon} \frac{\ip{\cT_N \bm{z}}{\xi_n}}{\sigma_n} \bm{v}_n}^2 = \sum_{\sigma_n > \epsilon} | \ip{\bm{z}}{\bm{v}_n}|^2 \leq \| \bm{z} \|^2.
}
Combining these two bounds now gives the first result.  For the second result, we first let $\bm{a}_N \in \bbC^N$ be the vector with $n^{\rth}$ entry $a_n = \ip{f}{\cS^{-1} \phi_n}$ for $n \in I_N$.  Then
\bes{
\| \bm{a} - \bm{a}^{\epsilon}_N \| \leq \sqrt{\sum_{n \in I \backslash I_N} | a_n |^2} + \| \bm{a}_N - \bm{x}^{\epsilon} \|.
} 
Hence it suffices to estimate $\| \bm{a}_N - \bm{x}^{\epsilon} \|$.  For this, we note that $f = \cS \cS^{-1} f = \cS_N \cS^{-1} f + (\cS - \cS_N) \cS^{-1} f$.  Since $\cS_N$ is self-adjoint and $\cS_{N} \xi_n = \cT_N \cT^*_N \cT_N \bm{v}_n = \sigma_n \cT_N \bm{v}_n = \sigma_n \xi_n$ we have
\bes{
\ip{f}{\xi_n} = \ip{\cS_N \cS^{-1} f}{\xi_n} + \ip{(\cS - \cS_N) \cS^{-1} f}{\xi_n} = \sigma_n \ip{\cS^{-1} f }{\xi_n} +  \ip{(\cS - \cS_N) \cS^{-1} f}{\xi_n} .
}
Therefore
\be{
\label{step3}
\bm{x}^{\epsilon} = \sum_{\sigma_n > \epsilon} \frac{ \ip{f}{\xi_n}}{\sigma_n} \bm{v}_n =  \sum_{\sigma_n > \epsilon} \ip{\cS^{-1} f}{\xi_n} \bm{v}_n + \sum_{\sigma_n > \epsilon} \frac{1}{\sigma_n}  \ip{(\cS - \cS_N) \cS^{-1} f}{\xi_n} \bm{v}_n.
}
Conversely, since $\bm{a}_N = \cT^*_N \cS^{-1} f$ we have $\ip{\bm{a}_N}{\bm{v}_n}= \ip{\cS^{-1} f}{\cT_N \bm{v}_n} = \ip{\cS^{-1} f }{\xi_n}$.  Hence
\be{
\label{step4}
\bm{a}_N = \sum_{n \in I_N} \ip{\bm{a}_N}{\bm{v}_n} \bm{v}_n = \sum_{n \in I_N} \ip{\cS^{-1} f}{\xi_n} \bm{v}_n.
}
Combining \R{step3} and \R{step4} now gives
\be{
\label{step5}
\| \bm{a}_N - \bm{x}^{\epsilon} \| \leq  \nm{\sum_{\sigma_n \leq \epsilon} \ip{\cS^{-1} f}{\xi_n} \bm{v}_n} + \nm{\sum_{\sigma_n > \epsilon} \frac{1}{\sigma_n}  \ip{(\cS - \cS_N) \cS^{-1} f}{\xi_n} \bm{v}_n}.
}
Consider the first term.  By orthogonality
\be{
\nm{\sum_{ \sigma_n \leq \epsilon} \ip{\cS^{-1} f}{\xi_n} \bm{v}_n}^2  \leq \epsilon \sum_{\sigma_n \leq \epsilon} \frac{1}{\sigma_n} | \ip{\cS^{-1} f}{\xi_n} |^2 \leq \epsilon \| \cS^{-1} f \|^2 \leq \epsilon/A \| \bm{a} \|^2. \label{bit1}
}
Here the final inequality follows from \R{framecond}.  Now consider the second term:
\eas{
 \nm{\sum_{\sigma_n > \epsilon} \frac{1}{\sigma_n}  \ip{(\cS - \cS_N) \cS^{-1} f}{\xi_n} \bm{v}_n}^2 &= \sum_{ \sigma_n > \epsilon} \frac{1}{\sigma^2_n} | \ip{(\cS - \cS_N) \cS^{-1} f}{\xi_n} |^2
 \\
 & \leq \frac{1}{\epsilon}  \sum_{ \sigma_n > \epsilon} \frac{1}{\sigma_n} | \ip{(\cS - \cS_N) \cS^{-1} f}{\xi_n} |^2
 \\
 & \leq \frac{1}{\epsilon} \| (\cS - \cS_N) \cS^{-1} f \|^2.
}
Observe that
\bes{
\nm{(\cS - \cS_N) \cS^{-1} f } = \nm{\sum_{n \in I \backslash I_N} a_n \phi_n} = \sup_{\substack{g \in \rH \\ g \neq 0}} \left \{ \frac{\left | \sum_{n > N} a_n \overline{\ip{g}{\phi_n}} \right | }{\| g \|} \right \} \leq \sqrt{B} \sqrt{\sum_{n \in I \backslash I_N} | a_n |^2},
}
and therefore
\bes{
 \nm{\sum_{\sigma_n > \epsilon} \frac{1}{\sigma_n}  \ip{(\cS - \cS_N) \cS^{-1} f}{\xi_n} \bm{v}_n}^2 \leq B/\epsilon \sum_{n \in I \backslash I_N} | a_n |^2.
}
Substituting this and \R{bit1} into \R{step5} gives the result.
}

Theorems \ref{t:Proj_err} and \ref{t:coefficients1} show a rather surprising conclusion.  Despite severe ill-conditioning of the Gram matrix, which led us to discard all of its singular values of size less than $\epsilon$, one still gets convergence of  $\cP^{\epsilon}_N f$ to within $\sqrt{\epsilon}$ of $f$.  Moreover, although the coefficients $\bm{x}^{\epsilon}$ may initially grow large (due to the $1/\sqrt{\epsilon}$ factor in \R{coeffics_bounded}), they too eventually converge to within $\sqrt{\epsilon}$ of the frame coefficients of $f$.

The underlying reason for this, as detailed in the following theorem, is that the absolute condition number of the mapping $\bm{y} = \{ \ip{f}{\phi_n} \}_{n \in I_N} \mapsto \cP^{\epsilon}_N f$ from the data $\bm{y}$ to the projection $\cP^{\epsilon}_N f$ is much smaller -- precisely, on the order of $1/\sqrt{\epsilon}$ as opposed to $1/\epsilon$ -- than the absolute condition number of the mapping from $\bm{y}$ to the coefficients $\bm{x}^{\epsilon}$ (which is just the minimal singular value of the SVD truncated Gram matrix $\bm{G}^{\epsilon}_N$; see \R{trunc_Gram_eps}).  Hence, whilst errors in the coefficients may be on the order of $1/\epsilon$, they result in much smaller errors in the projection $\cP^{\epsilon}_N$.

\thm{
\label{t:Proj_Condit}
The absolute condition number of the mapping $\bbC^N \rightarrow \rH^{\epsilon}_N$, $\bm{y} \mapsto \cT_N (\bm{G}^{\epsilon}_N)^{\dag} \bm{y}$ satisfies
$\kappa \leq \min \left \{ 1/\sqrt{\sigma_{\min}(\bm{G}_N)}, 1/\sqrt{\epsilon} \right \}$,
where $\sigma_{\min}(\bm{G}_N)$ is the minimal singular value of the Gram matrix $\bm{G}_N$.
}
\prf{
By linearity, the condition number of the mapping $\bm{y} \mapsto \cT_N (\bm{G}^{\epsilon}_N)^{\dag} \bm{y}$ is 
\bes{
\kappa = \max_{\substack{\bm{y} \in \bbC^N \\ \| \bm{y} \|=1}} \| \cT_N (\bm{G}^{\epsilon}_N)^{\dag} \bm{y} \|.
}
We have
\bes{
 \| \cT_N (\bm{G}^{\epsilon}_N)^{\dag} \bm{y} \|^2 = \ip{ \bm{y}}{(\bm{G}^{\epsilon}_N)^{\dag}\bm{G}_N (\bm{G}^{\epsilon}_N)^{\dag} \bm{y}} = \ip{\bm{y}}{(\bm{G}^{\epsilon}_N)^{\dag} \bm{y}} = \sum_{\sigma_n > \epsilon} \frac{|\ip{\bm{y}}{\bm{v}_n}|^2}{\sigma_n}.
}
This gives $\| \cT_N (\bm{G}^{\epsilon}_N)^{\dag} \bm{y} \|^2 \leq 1/\min \{ \sigma_n : \sigma_n > \epsilon \} \| \bm{y} \|^2$, and the result follows.
}

\rem{
Theorems \ref{t:Proj_err} and \ref{t:coefficients1} assert convergence to within $\sqrt{\epsilon}$ of $f$ only.  Using spectral theory techniques it can be shown that $\bm{a}^{\epsilon}_N \rightarrow \bm{a}$ as $N \rightarrow \infty$ \cite[Thm.\ 5.17]{HarrisonThesis}, i.e.\ the regularized coefficients converge to the frame coefficients in the canonical dual frame.  Since $\| f - \cP^{\epsilon}_{N} f \| = \|\cT (\bm{a} -\bm{a}^{\epsilon}_N) \| \leq \sqrt{B} \| \bm{a} - \bm{a}^{\epsilon}_N \|$ this also gives $\cP^{\epsilon}_{N} f \rightarrow f$, albeit at a rate that can be arbitrarily slow.  Of course, in finite precision calculations convergence beyond $\ord{\sqrt{\epsilon}}$ will not be expected, due to the condition number of the mapping (Theorem \ref{t:Proj_Condit}).  We note also that convergence down to $\ord{\sqrt{\epsilon}}$ has previously been observed empirically in \cite{HarrisonThesis,StrohmerNonuniform}.  The main contribution of Theorems \ref{t:Proj_err} and \ref{t:coefficients1} is that they provide explicit bounds which can be used to estimate the rate of decay of these errors in the regime where they are larger than $\sqrt{\epsilon}$. 
}

\subsection{Discussion and examples}
\label{ss:example_convergence}

We now consider Theorems \ref{t:Proj_err}--\ref{t:Proj_Condit} in relation to the examples of \S \ref{sect:examples}.  We focus on two issues: (i) the convergence of the projection $\cP^{\epsilon}_N f$, and (ii) the behaviour of the coefficients $\bm{x}^{\epsilon}$.

For proofs of the various results presented, we refer to the supplementary material.

First, notice that for small $N$ -- specifically, for $N$ such that $\sigma_{\min}(\bm{G}_N) \geq \epsilon$ -- we have $\cP^{\epsilon}_N = \cP_N$.  Hence, the truncated SVD projection initially behaves like the exact projection $\cP_N$.  However, beyond this point the convergence begins to differ.  If $\bm{x} \in \bbC^N$ are the coefficients of $\cP_N f$ then setting $\bm{z} = \bm{x}$ in \R{Proj_err} gives
\bes{
\| f - \cP^{\epsilon}_N f \| \leq \| f - \cP_N f \| + \sqrt{\epsilon} \| \bm{x} \|.
}
As discussed in \S \ref{ss:impossible}, the term $\|\bm{ x} \|$ often grows rapidly in $N$.  Hence as $N$ increases, the right-hand side of the above inequality may begin to diverge.

However, since $\Phi$ is a frame there are infinitely-many sequences  of coefficients $\bm{c} \in \ell^2(I)$ such that $f = \cT \bm{c}$.  Suppose that $\bm{z} = \{ c_n \}_{n \in I_N}$.  Then Theorem \ref{t:Proj_err} gives
\be{
\label{bounded_err}
\| f - \cP^{\epsilon}_N f \| \leq \nm{ f - \cT_N \bm{c} } + \sqrt{\epsilon} \| \bm{c} \| = \nm{  \sum_{n \in I \backslash I_N} c_n \phi_n }  + \sqrt{\epsilon} \| \bm{c} \| .
}
The term $\sqrt{\epsilon} \| \bm{c} \|$ is independent of $N$, while the other term tends to zero as $N \rightarrow \infty$.  Hence, the rate of decay of the error down to $\sqrt{\epsilon}$ is bounded by how well $f$ can be expanded in the frame using coefficients $\bm{c}$ with small norms.  In the examples below, we show that there always exist coefficient sequences $\bm{c}$ that achieve favourable rates of decay.  Therefore, while the truncated SVD projection $\cP^{\epsilon}_N f$ may not achieve the same error decay as the exact projection $\cP_N f$, we can often expect good accuracy.

\vspace{1pc}\noindent
\textbf{Fourier frames for complex geometries.}  Consider the frame of Example 1.  We have the following:

\prop{
\label{l:FE_bounded_conv}
Let $\Omega \subseteq (-1,1)^d$ be a compact, Lipschitz domain and consider the frame \R{FE_frame}.   If $f \in \rH^{kd}(\Omega)$ then there exists a $\bm{c} \in \ell^2(I)$ such that 
\bes{
\| f - \cT_N \bm{c} \| \leq C N^{-k} \| f \|_{\rH^{kd}(-1,1)^d},\qquad \| \bm{c} \| \leq C \| f \|_{\rH^{kd}(-1,1)^d},
}
where $C > 0$ is independent of $f$ and $N$.  In particular, for the exact projection
\bes{
\| f - \cP_{N} f \| \leq C N^{-k} \| f \|_{\rH^{kd}(-1,1)^d},
}
and for the regularized projection,
\be{
\label{FE_num_conv}
\| f - \cP^{\epsilon}_N f \| \leq C \left ( N^{-k} + \sqrt{\epsilon} \right ) \| f \|_{\rH^k(-1,1)^d}.
}
}

\begin{figure}
\begin{center}
$\begin{array}{ccc}
\includegraphics[width=6.0cm]{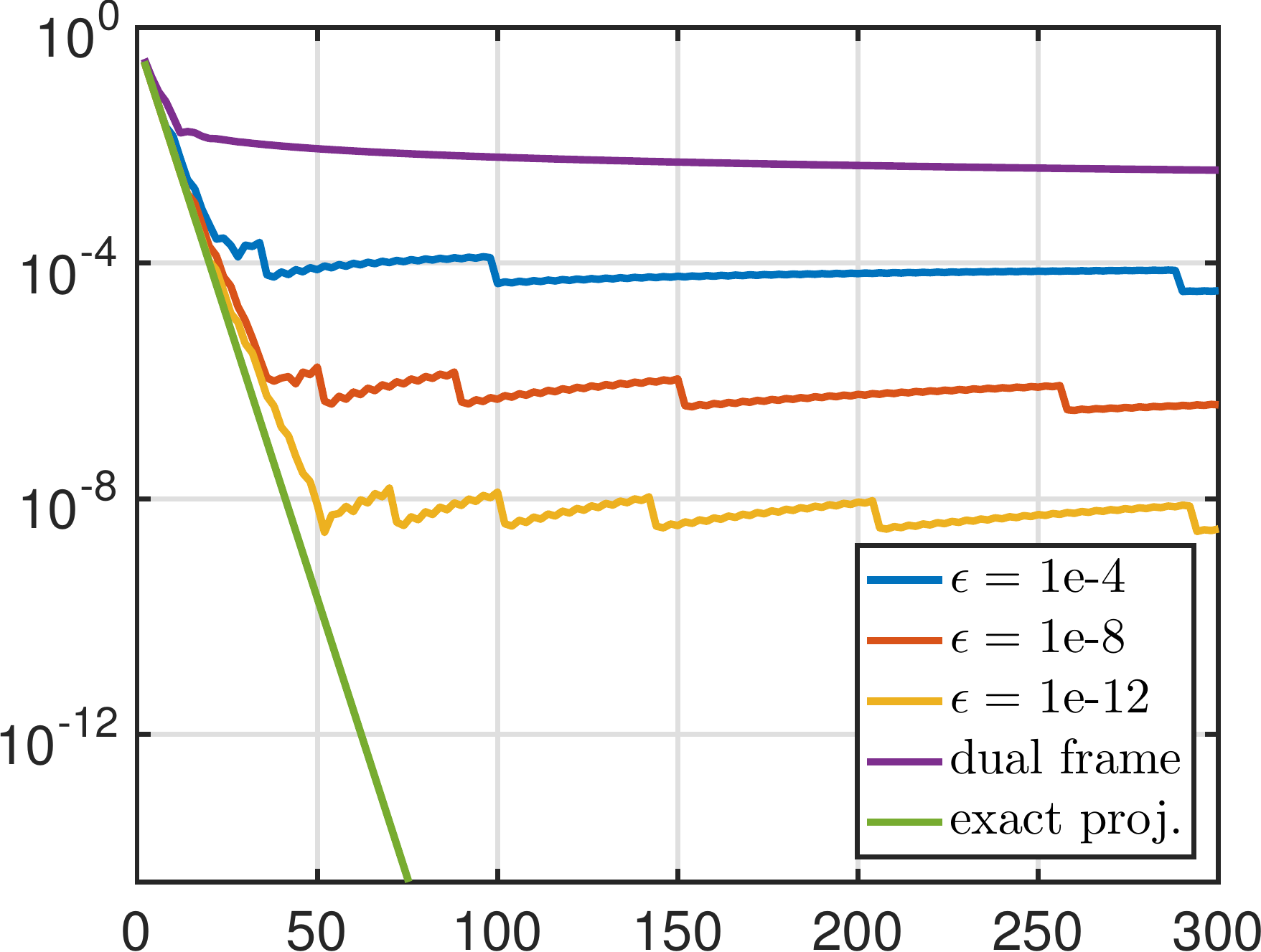} && \includegraphics[width=6.0cm]{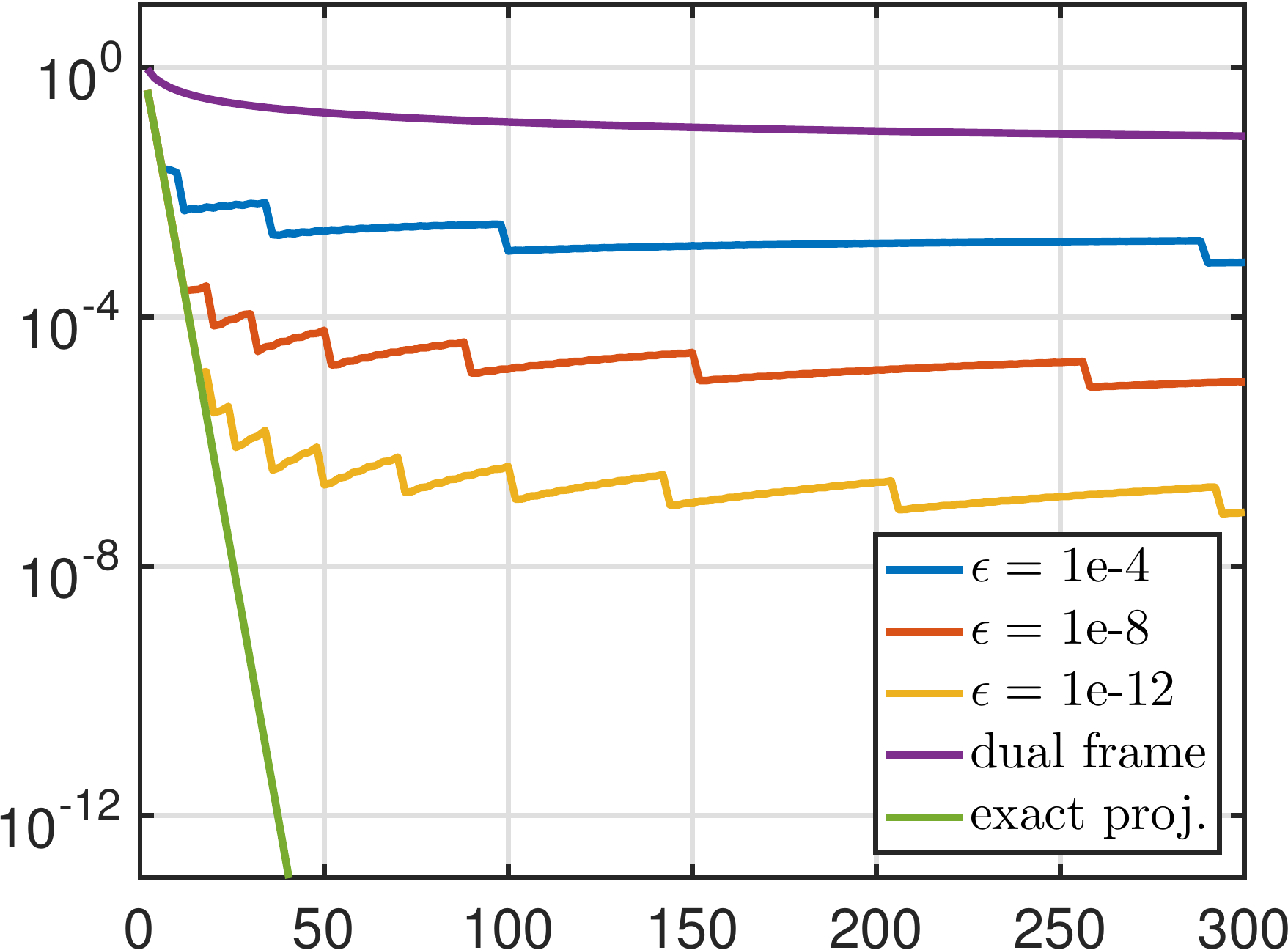} \\
\includegraphics[width=6.0cm]{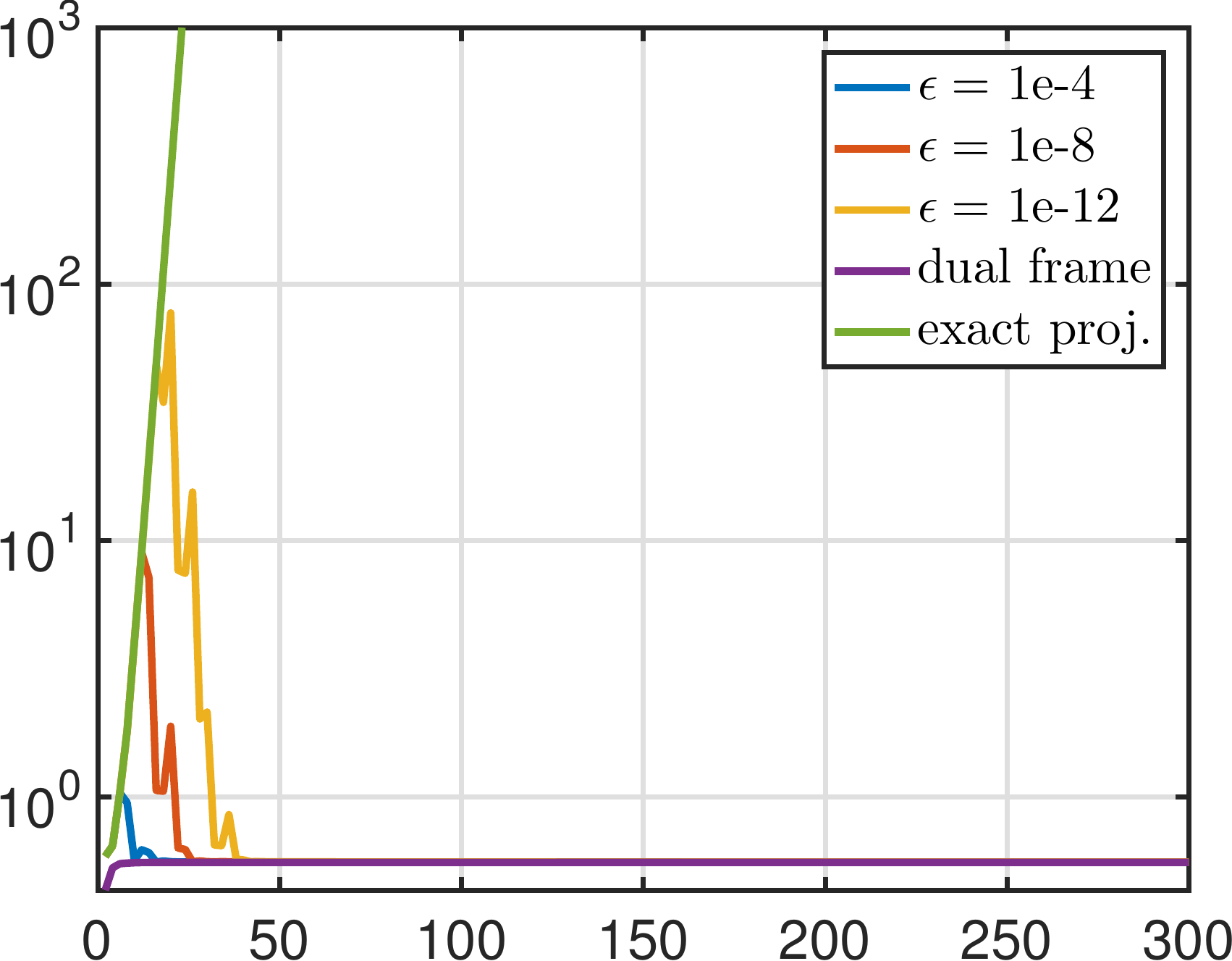} &&\includegraphics[width=6.0cm]{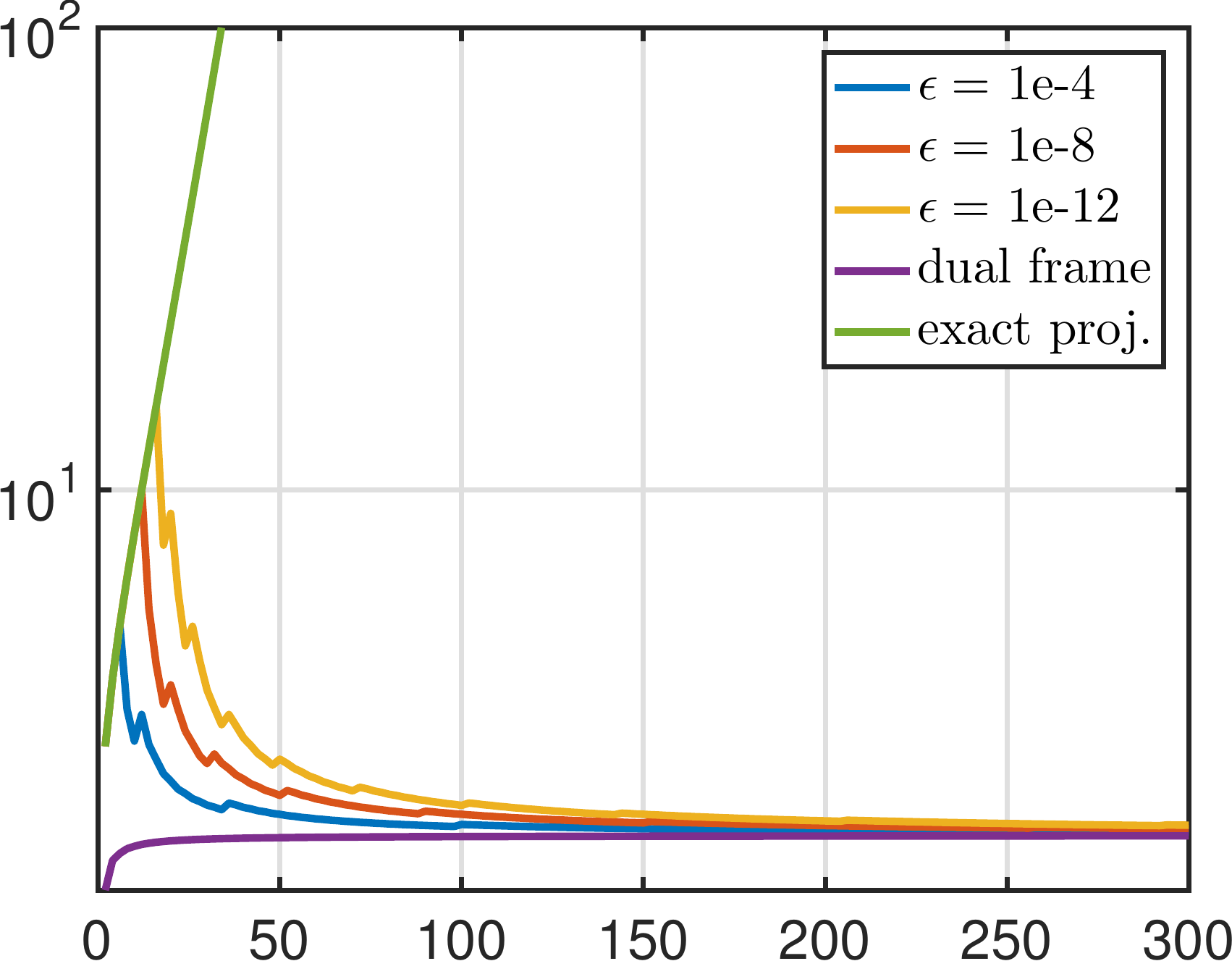} \\
f_1(t) = \frac{1}{1+25 t^2} && f_2(t) = \frac{1}{3/4 - t} 
\end{array}$
\caption{Comparison of the exact projection $\cP_N f$, the truncated SVD projections $\cP^{\epsilon}_N f$ and the canonical dual frame expansion $\sum_{n \in I_N} a_n \phi_n$ for Example 1 with $T = 2$.  Top row: $\rL^2$-norm errors versus $N$.  Bottom row: the norms of the coefficient vectors versus $N$.  The exact projection was computed using additional precision.  All other approximations were computed in double precision.}
\label{f:FE_Err_Coeff}
\end{center}
\end{figure}

This proposition asserts that there are bounded coefficient vectors for which the error $\| f - \cT_N \bm{c} \|$ decays at an arbitrarily-high algebraic rate.  In particular, the error of the exact projection decays spectrally fast in $N$.  On the other hand, for the regularized projection we have the following.  First, for smooth functions $f$, the error decays rapidly in $N$ when  $\| f - \cP^{\epsilon}_N f \| \gg \sqrt{\epsilon}$.  Second, if the $k^{\rth}$ derivatives of $f$ grow rapidly in $k$ then the rate of decrease of the error may lessen as $\| f - \cP^{\epsilon}_N f \|$ approaches $\sqrt{\epsilon}$, due to the presence of the Sobolev norm in the error bound \R{FE_num_conv}.
This effect is shown in Fig.\ \ref{f:FE_Err_Coeff}.  Note that the derivatives of the function $f_1$ grow slowly with $k$, whereas the derivatives of $f_2$ grow much more rapidly.  As predicted by \R{FE_num_conv}, there is substantially less effect in replacing $\cP_N$ by the regularized projection $\cP^{\epsilon}_N$ in the case of $f_1$ than in the case of $f_2$. This aside, Fig.\ \ref{f:FE_Err_Coeff} also shows the slow convergence of the canonical dual frame expansion \R{dual_rep}, thus confirming the discussion in \S \ref{sect:examples}, and the norms of the various coefficient vectors.  These results are in good agreement with Theorem \ref{t:coefficients1}: the coefficients of $\cP^{\epsilon}_N f$ initially grow large, but after a certain point they begin to decay to the limiting value $\nm{\bm{a}}$.

\begin{figure}
\begin{center}
$\begin{array}{ccc}
\includegraphics[width=6.0cm]{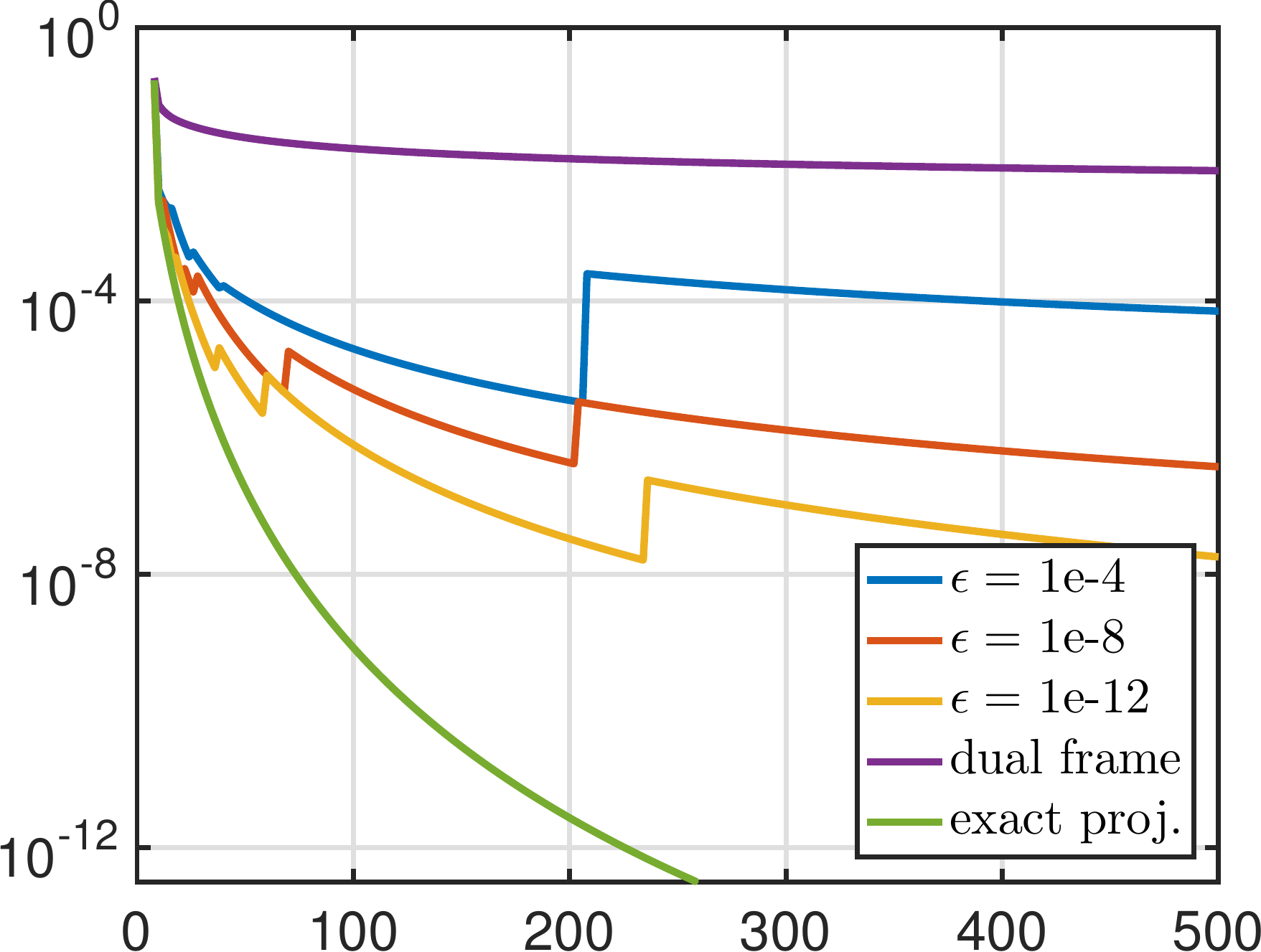} &&  \includegraphics[width=6.0cm]{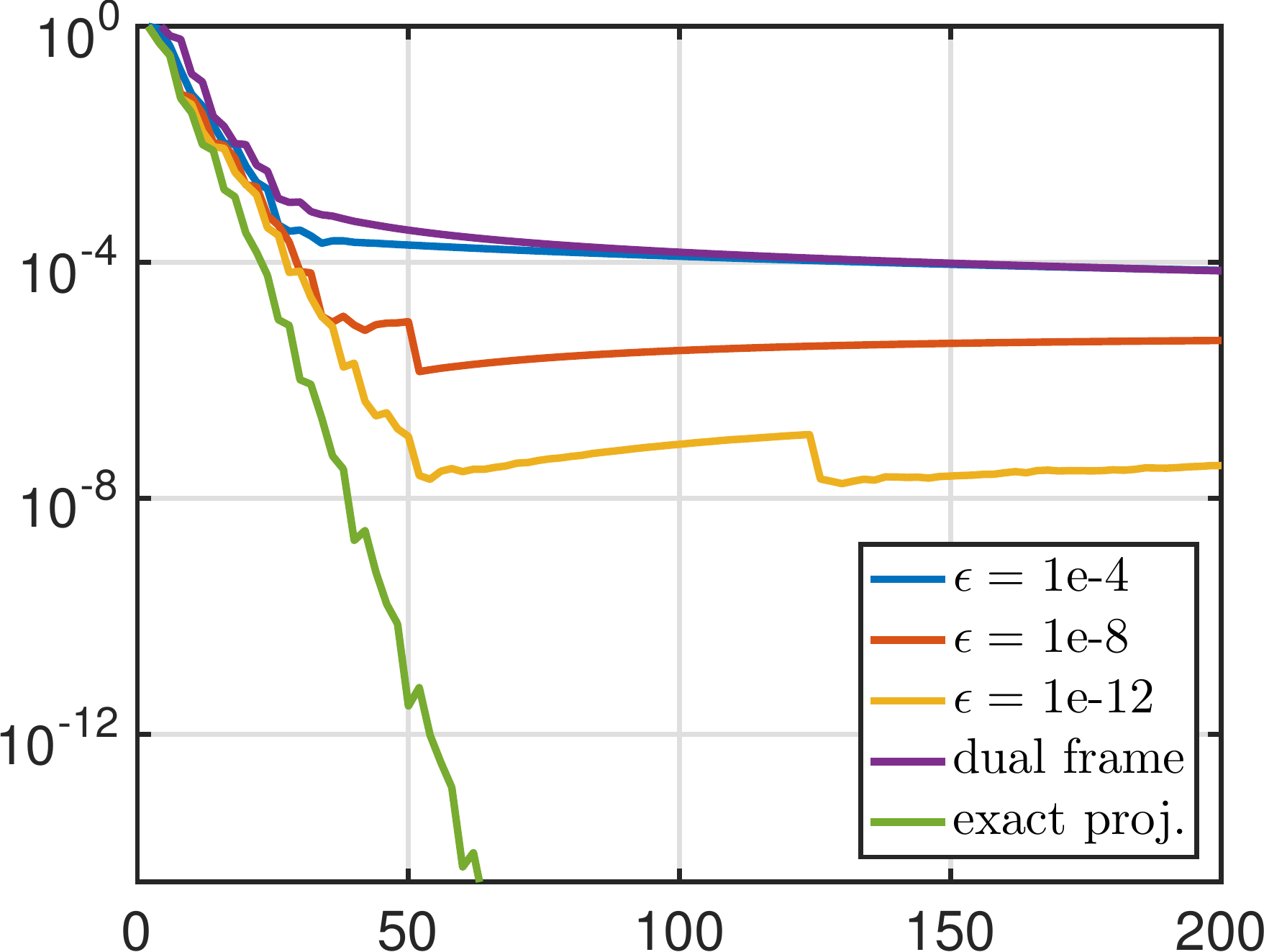}\\
f_1(t) = \frac{1}{10-9 t} && f_2(t) = \E^{\sin(3t+1/2)}\sqrt{1+t} +\cos(5t)
\end{array}$
\caption{$\rL^2$-norm errors versus $N$ for the exact projection $\cP_N f$, the truncated SVD projections $\cP^{\epsilon}_N f$ and the canonical dual frame expansion $\sum_{n \in I_N} a_n \phi_n$ for Example 2 with $K = 8$ (left) and Example 3 with $\alpha = 1/2$ (right).  The exact projection was computed using additional precision.  All other approximations were computed in double precision.
%Top row: the errors versus $N$.  Bottom row: the norms of the coefficient vectors versus $N$.
}
\label{f:FP_Err_Coeff}
\end{center}
\end{figure}

\vspace{1pc}\noindent
\textbf{Augmented Fourier basis.} Consider the frame of Example 2.
In this case we observe algebraic convergence at a rate determined by the number of polynomials $K$ added to the Fourier basis.  See also Fig.\ \ref{f:FP_Err_Coeff}. 
\prop{
\label{l:AugFourier_conv}
Let $K \in \bbN$ be fixed and consider the frame \R{AugFourier_frame}.  If $f \in \rH^k(-1,1)$ for $0 \leq k \leq K$ then there exists a $\bm{c} \in \ell^2(I)$ such that
\bes{
\| f - \cT_N \bm{c} \| \leq C N^{-k} \| f \|_{\rH^k(-1,1)},\qquad \| \bm{c} \| \leq C \| f \|_{\rH^k(-1,1)},
}
where $C > 0$ is independent of $f$ and $N$.  In particular,
\eas{
\| f - \cP_N f \| &\leq C N^{-k}  \| f \|_{\rH^k(-1,1)},
\\
\| f - \cP^{\epsilon}_N f \| &\leq C \left ( N^{-k} + \sqrt{\epsilon} \right ) \| f \|_{\rH^k(-1,1)}.
}
}

%This result establishes algebraic convergence in this frame.  

\vspace{1pc}\noindent
\textbf{Polynomials plus modified polynomials.} For the frame of Example 3 we have:
\prop{
\label{l:PMP_conv}
Let $\cQ_{N} : \rL^2(-1,1) \rightarrow \rL^2(-1,1)$ be the orthogonal projection onto $\bbP_{N/2-1}$.  If $f(t) = w(t) g(t) + h(t)$, then there exists a $\bm{c} \in \ell^2(I)$ such that
\bes{
\| f - \cT_N \bm{c} \| \leq w_{\max} \| g - \cQ_{N} g \| + \| h - \cQ_{N} h\|,\qquad \| \bm{c} \| \leq  \| g \| + \| h \| ,
}
where $w_{\max} = \esssup_{t \in (-1,1)} | w(t) |$.  In particular,
\eas{
\| f - \cP_N f \| &\leq w_{\max} \| g - \cQ_{N} g \| + \| h - \cQ_{N} h\|,
\\
\| f - \cP^{\epsilon}_N f \| &\leq w_{\max} \| g - \cQ_{N} g \| + \| h - \cQ_{N} h\| + \sqrt{\epsilon} \left ( \| g \| + \| h \| \right ).
}
}

This result implies the convergence of the regularized projection $\cP^{\epsilon}_{N} f$ is spectral in the factors $g$ and $h$.  In particular, if $g$ and $h$ are infinitely differentiable then one has superalgebraic convergence down to $\sqrt{\epsilon}$, and if $g$ and $h$ are analytic, one has geometric convergence down to $\sqrt{\epsilon}$.  See Fig.\ \ref{f:FP_Err_Coeff}.  This figure suggests geometric convergence, in agreement with Proposition \ref{l:PMP_conv}, but with a somewhat reduced exponent over that of the exact projection $\cP_N f$.  In other words, there exist coefficient vectors in the frame which yield faster geometric convergence, but are too large to be obtained as solutions of the regularized system.

\section{Well-conditioned frame approximations}\label{s:numstableapprox}

As shown in the previous section, the accuracy of the numerical frame projection $\cP^{\epsilon}_{N}$ is limited to  $\ord{\sqrt{\epsilon}}$, and the relevant condition number is at most $1/\sqrt{\epsilon}$.  We close this paper with a brief description of an approximation that achieves $\ord{\epsilon}$ accuracy and has a bounded condition number independent of $\epsilon$.  This approximation is based on oversampling.  Specifically, rather than solving the square system \R{GramLinSyst}, we consider the $M \times N$ system
\be{
\label{OSGram}
\bm{G}_{M,N} \bm{x} \approx \bm{y},\qquad \bm{y} = \{ \ip{f}{\phi_n} \}_{n \in I_M},
}
where $\bm{G}_{M,N} = \{ \ip{\phi_n}{\phi_m} \}_{m \in I_M,n \in I_N} \in \bbC^{M \times N}$.  Note that where $\bm{G}_N$ corresponds to the finite section of the Gram operator $\cG$, $\bm{G}_{M,N}$ corresponds to a so-called \textit{uneven section}.  Uneven sections are known to be useful alternatives to finite sections in computational spectral theory \cite{bottcher1996,strohmer,hansen2008,hansen2011,lindner2006} and, more recently, sampling theory \cite{BAACHShannon,BAACHGSCS,BAACHOptimality}.  Much the same is true in this instance.

Since Gram operators of frames are singular, $\bm{G}_{M,N}$ remains ill-conditioned even when $M \geq N$.  Hence we consider the regularized projection
\bes{
\bm{x}^{\epsilon} = (\bm{G}^{\epsilon}_{M,N})^{\dag} \bm{y},\qquad \cP^{\epsilon}_{M,N} f = \sum_{n \in I_N} (\bm{x}^{\epsilon})_n \phi_n.
}
where $\bm{G}^{\epsilon}_{M,N}$ is obtained by discarding all its singular values of $\bm{G}_{M,N}$ below $\epsilon$. 

Given a sufficient amount of oversampling, this approximation is well conditioned and achieves $\ord{\epsilon}$ accuracy.  Specifically, it can be shown that the condition number of the mapping $\bm{y} \mapsto \cT_N (\bm{G}^{\epsilon}_{M,N})^{\dag} \bm{y}$ is bounded by $C^{\epsilon}_{M,N} / \sqrt{A}$, and approximation $\cP^{\epsilon}_{M,N} f$ satisfies the error bound 
\be{
\label{OSerr}
\| f - \cP^{\epsilon}_{M,N} f \| \leq \left ( 1 + \sqrt{B/A} C^{\epsilon}_{M,N} \right ) \inf \left \{  \| f - \cT_{N} \bm{z} \| + \epsilon / \sqrt{B} \| \bm{z} \| : \bm{z} \in \bbC^N \right \}.
}
Here $C^{\epsilon}_{M,N}$ is a constant depending on $\Phi$, $\epsilon$, $M$ and $N$.  It satisfies
\bes{
\limsup_{M \rightarrow \infty} C^{\epsilon}_{M,N} \leq 1,
}
for any fixed $N$ and $\epsilon$.  In other words, this constant can be made arbitrarily close to one with sufficiently large $M$.  Note that the inequality \R{OSerr} is similar to the error bound \R{Proj_err} for the projection $\cP^{\epsilon}_{N} f$, except for appearance of $\epsilon$ in place of $\sqrt{\epsilon}$.  Hence the error now decays down to $\ord{\epsilon}$, with, as before, the rate of decay of the error being dictated by existence of expansions in the frame with small-norm coefficients.

The proofs of these results are given in a companion paper \cite{BADHFramesPart2}, which develops on this paper by introducing a general framework for frame approximations from `indirect' data.  This framework also includes discrete function samples, for example, which are typically much more convenient to work with than \R{OSGram} since they do not require evaluations of inner products.  

\begin{figure}
\begin{center}
$\begin{array}{ccc}
\includegraphics[width=5.0cm]{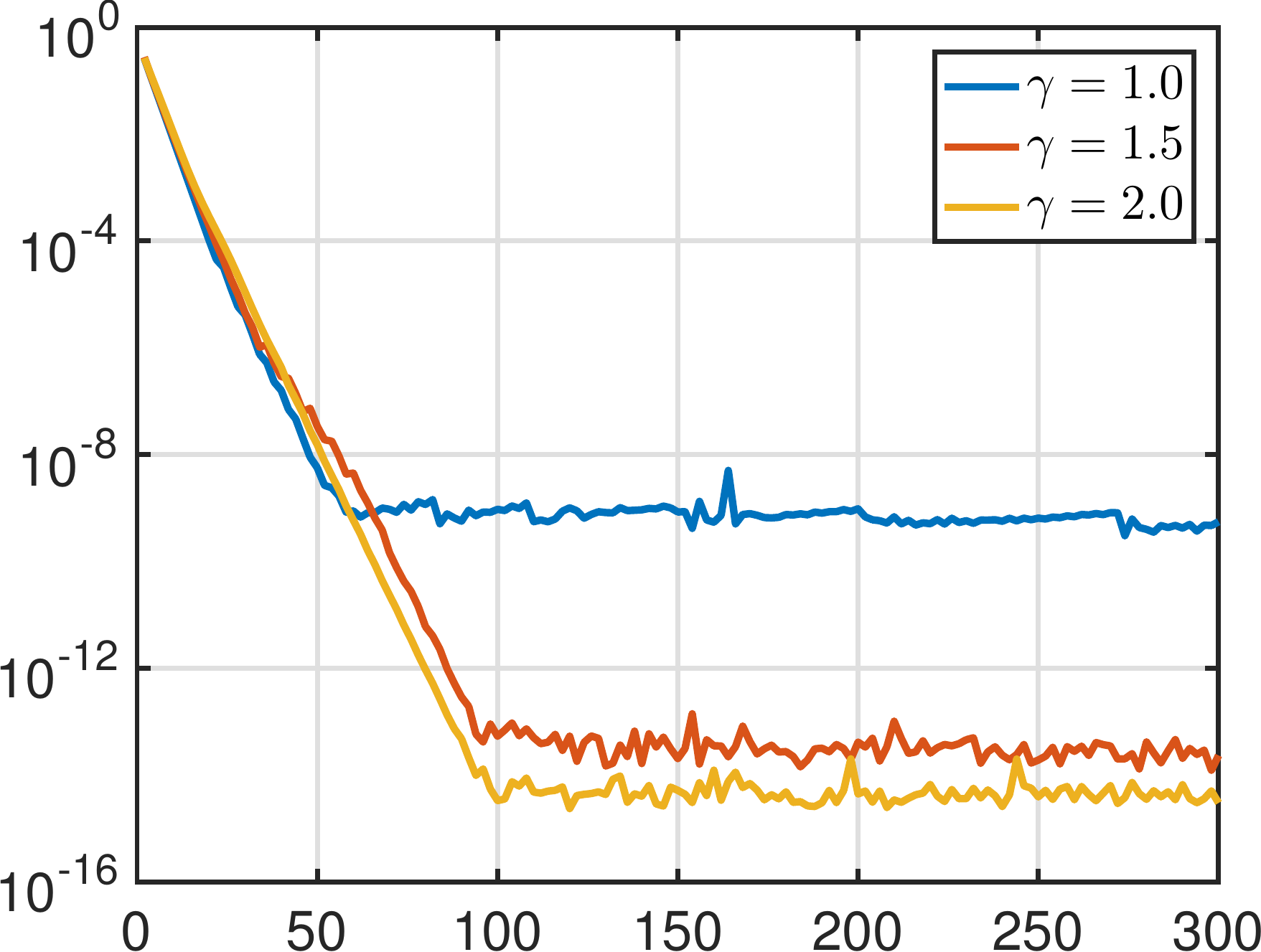} &
\includegraphics[width=5.0cm]{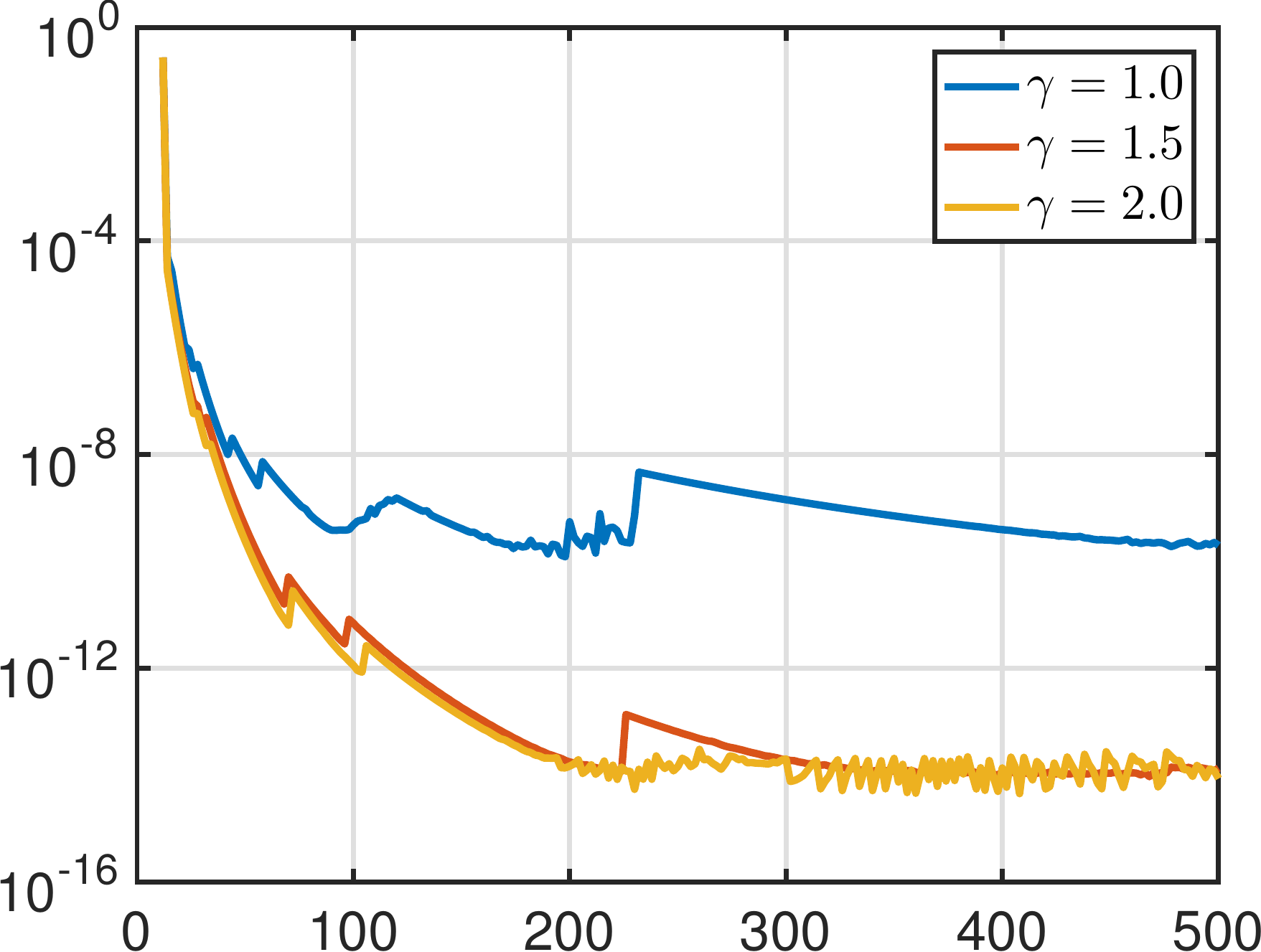}&
\includegraphics[width=5.0cm]{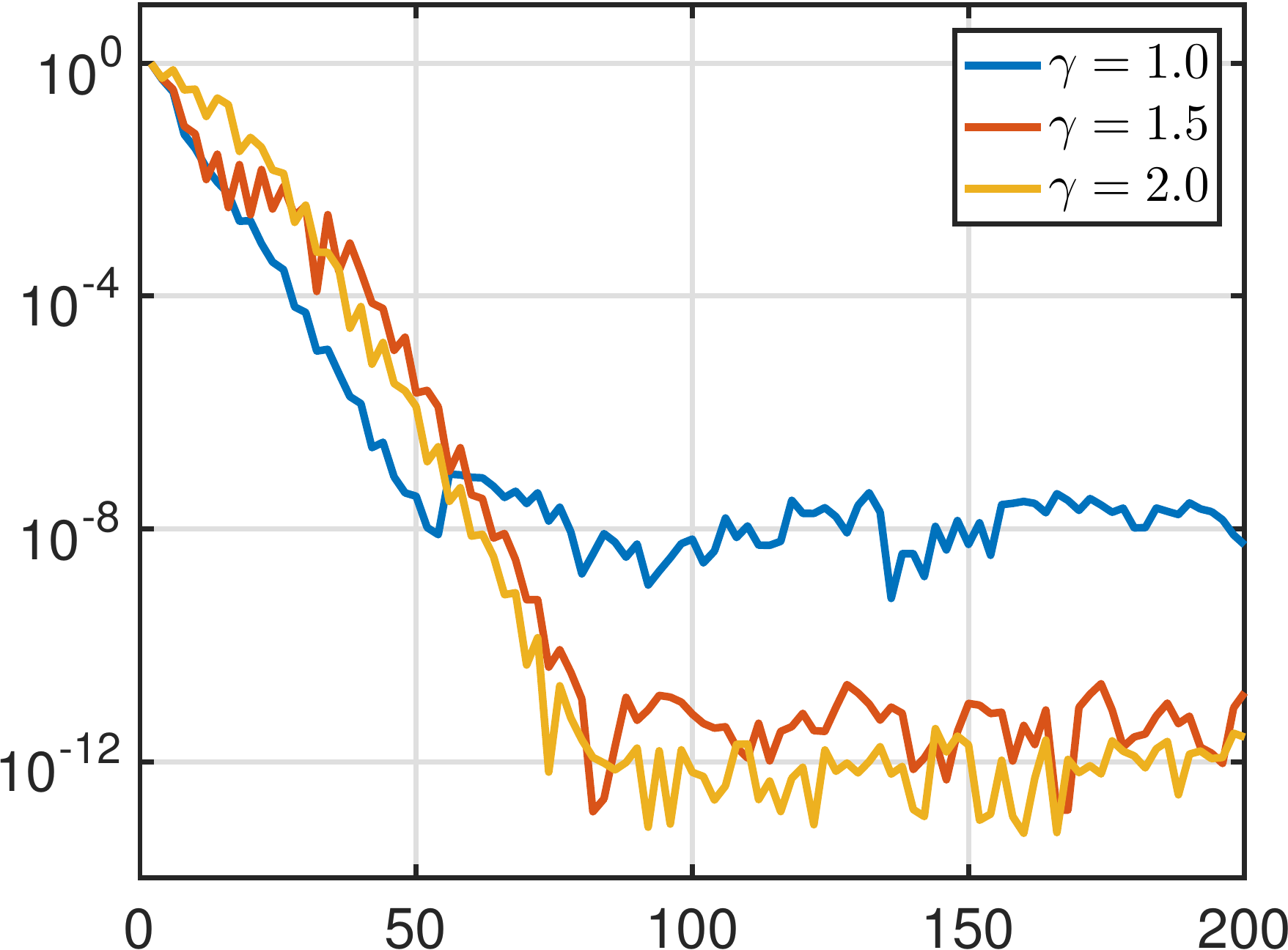}
\\
\includegraphics[width=5.0cm]{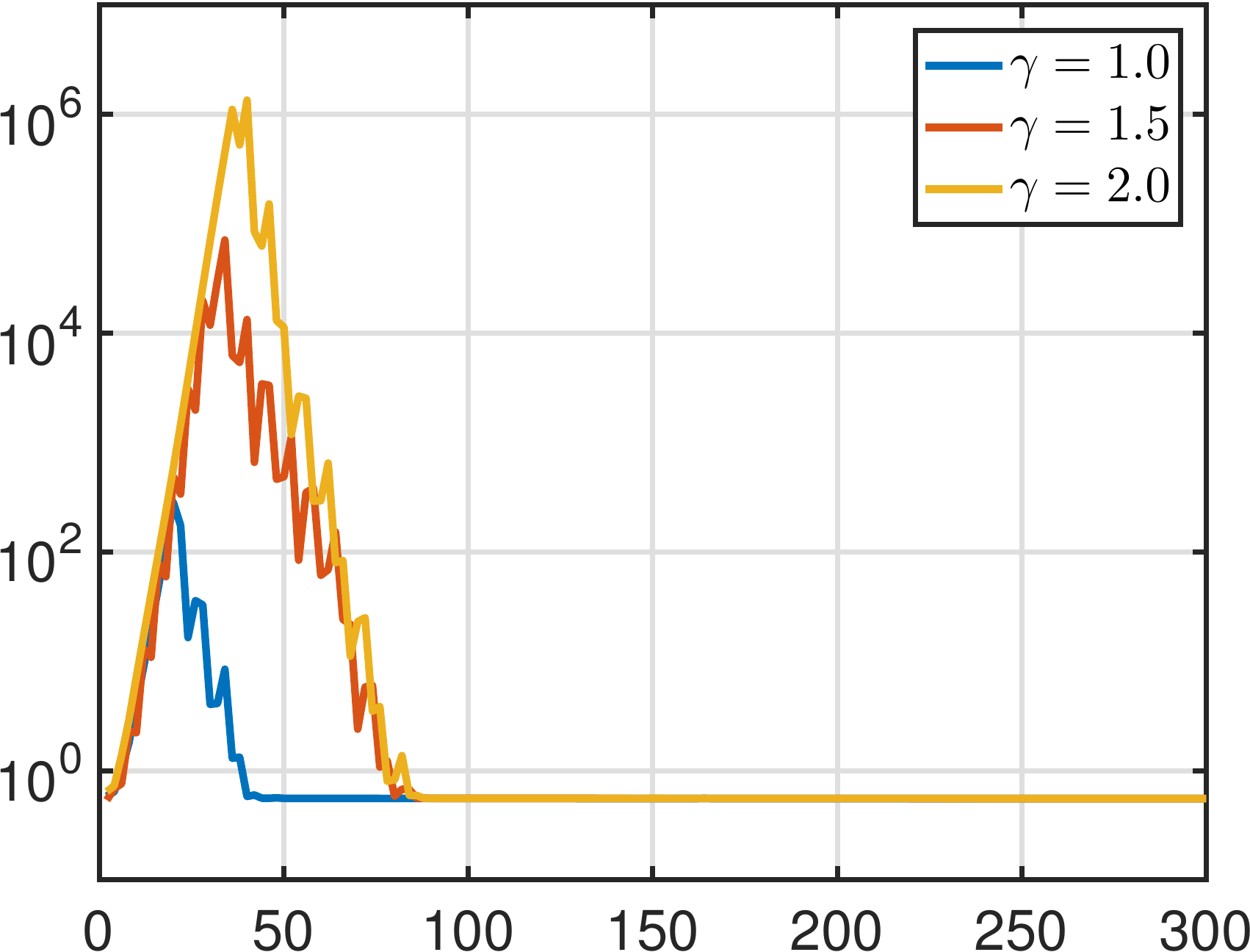}&
\includegraphics[width=5.0cm]{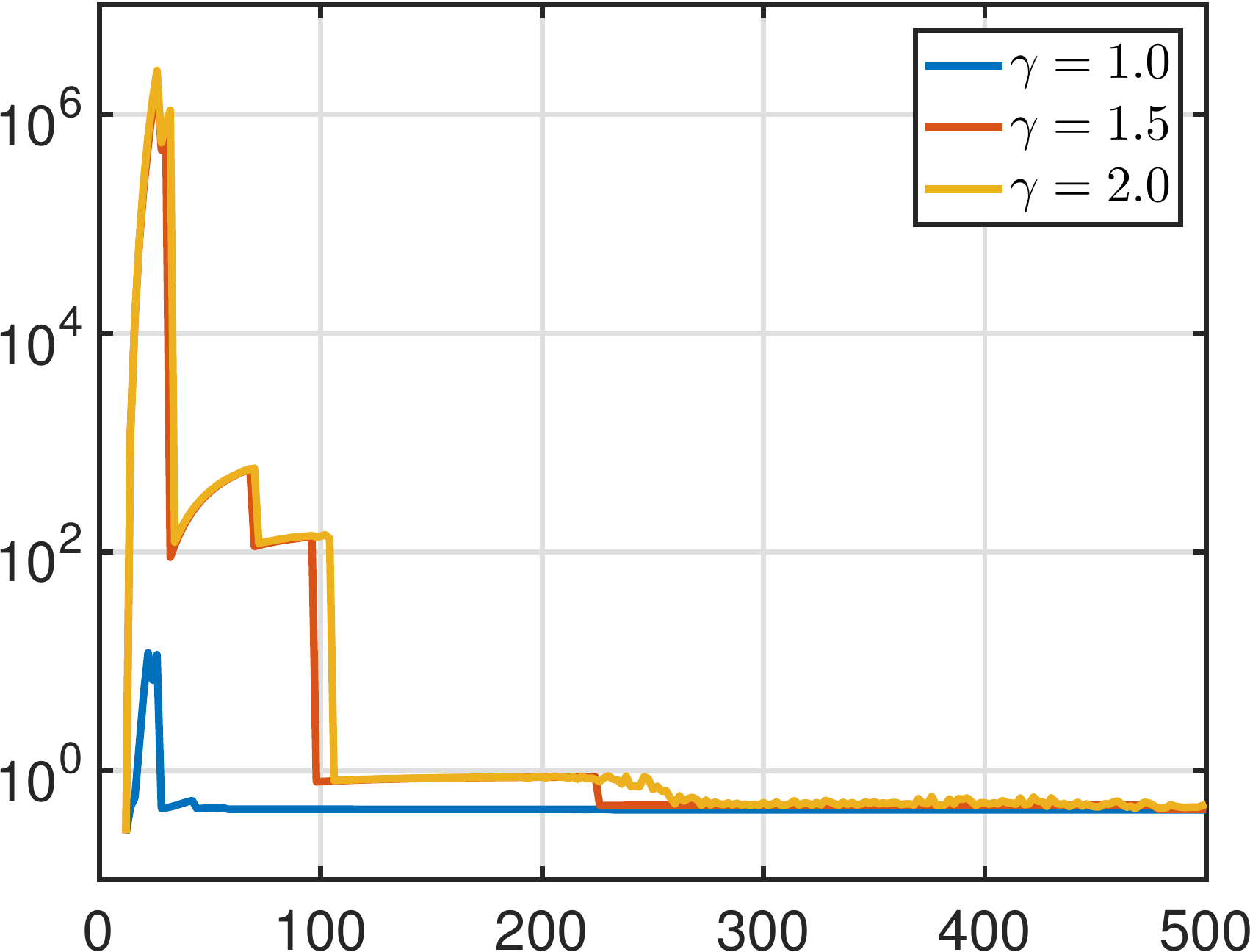}&
\includegraphics[width=5.0cm]{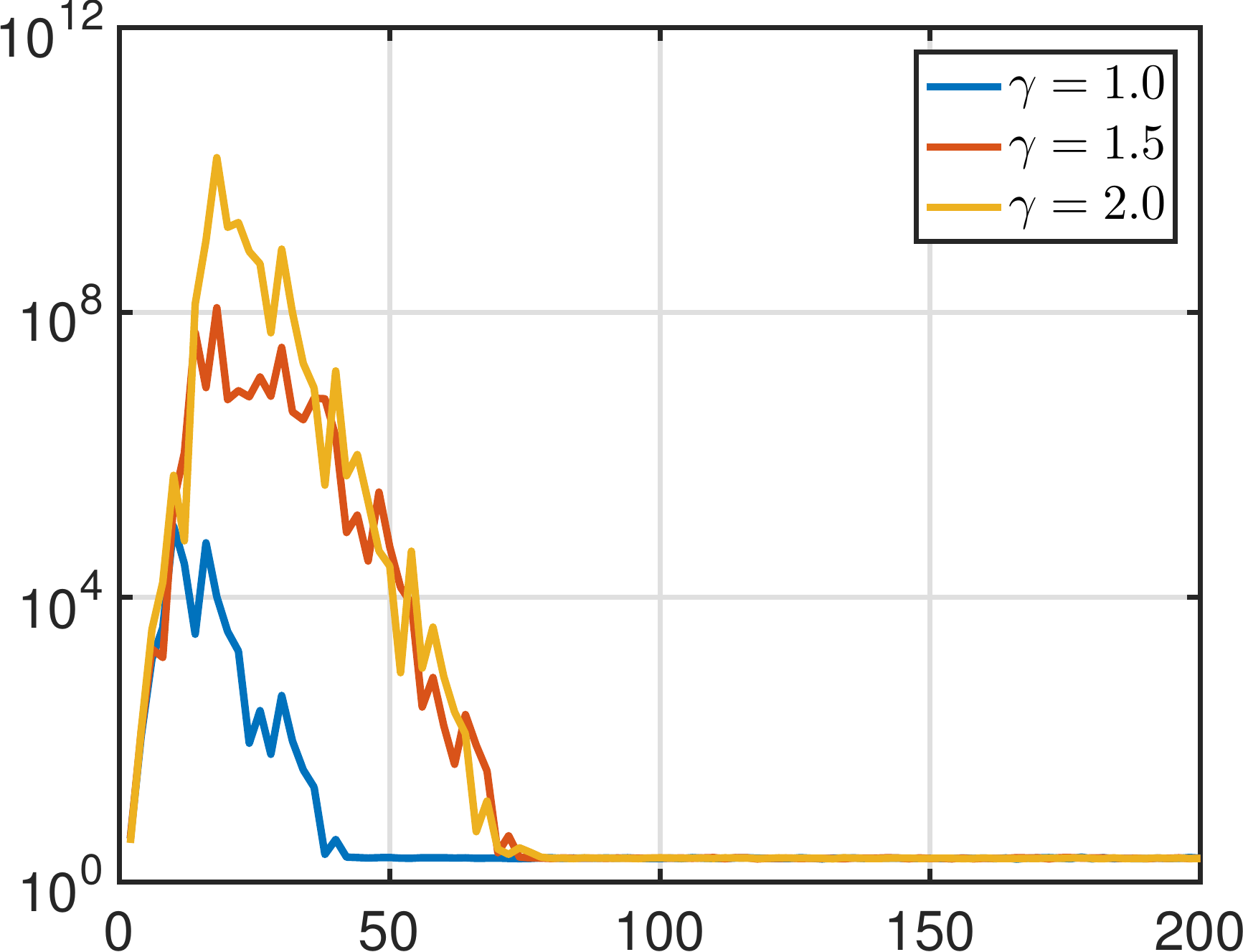}  \\
\mbox{{\small Example 1, $T = 2$}} & \mbox{{\small Example 2, $k = 12$}} & \mbox{{\small Example 3, $\alpha = 1/2$}} \\
\end{array}$
\caption{
Comparison of the projection $\cP^{\epsilon}_{\gamma N , N} f$ for different values of the parameter $\gamma$.  Top row: $\rL^2$-norm errors versus $N$.  Bottom row: the norms of the coefficient vectors versus $N$.  The functions used were  $f(t) = \frac{1}{1+25 t^2} $ (left), $ f(t) = \frac{1}{5 - 4 t}$ (middle) and $f(t) = \E^{\sin(3t+1/2)}\sqrt{1+t} +\cos(5t)$ (right).
}
\label{f:least_squares}
\end{center}
\end{figure}

In Fig.\ \ref{f:least_squares} we present several numerical results for the Examples 1--3 based on \R{OSGram}.  These results illustrate that with a mild amount oversampling one can obtain a much more accurate numerical frame approximation than when $M = N$, which is the case considered previously (see \cite{FEParameterSelection} for some further analysis in the case of Example 1).  Similar results can also be obtained via least-squares fitting with discrete function samples taken, for example, on an equally-spaced grid of $M \approx 2 N$ points \cite{BADHFramesPart2}. A fast algorithm to do so is given in \cite{matthysen2017fastfe2d}, and an illustration is shown in Fig.\ \ref{f:northpole}.

\section{Conclusions and further research}\label{s:conclusions}

The concern of this paper has been computing numerical approximations using frames, in particular, orthogonal projections in the span of the first $N$ frame elements.  There are four main conclusions.  First, the Gram matrices of truncated frames are necessarily ill-conditioned.  As shown, such ill-conditioning is due to approximation of the singular Gram operator by its finite section, and can be arbitrarily bad depending on the frame.  Second, the orthogonal projection typically cannot be computed, since its coefficients generally grow rapidly in $N$.  However, using regularization it is possible to find an approximation in a truncated frame whose coefficients are bounded (by at most $\ord{1/\sqrt{\epsilon}}$ and eventually are $\ord{1}$), has condition number at most $1/\sqrt{\epsilon}$, and which converges down to an error of $\ord{\sqrt{\epsilon}}$.  The convergence of this approximation depends on how well $f$ can be approximated by finite expansions in the frame with small-norm coefficients.  Fourth, and finally, using oversampling one can construct an approximation that has $\ord{1}$ condition number and is convergent down to $\ord{\epsilon}$.

The overall conclusion of this paper is that satisfactory approximations can be computed in certain finite systems with near-linear dependencies.  We stress that the frame condition is important in this regard.  The monomials $\Phi_N = \{ 1 ,x , x^2,\ldots x^{N-1} \}$ are nearly-linearly dependent for large $N$, but may not give good numerical approximations, since a smooth function does not necessarily have a representation in this system with small-norm coefficients.  On the other hand, a frame guarantees at least one approximation with small-norm coefficients -- namely, the truncated canonical dual frame expansion -- although, as seen, better approximations often exist.  At this stage we stress again that the approach developed in this paper is recommended only when the dual frame expansion gives unsatisfactory convergence, such as in Examples 1--3.  Otherwise, approximation via this expansion is arguably preferable (recall \S \ref{ss:background}).

There are a number of topics not considered in this paper.  First, our main results apply to general frames.  Yet there may be theoretical and practical advantages to considering classes of frames with particular properties; for instance, so-called \textit{localized frames} (see \cite{BalazsGrochenigLocalized,FornasierGrochenigLocal,GrochenigLocalization,SongGelbInverse} and references therein).

% is a more detailed analysis of the least-squares frame approximations, considered briefly in \S \ref{s:numstableapprox}.  This also includes the case of approximations from `indirect' data, e.g.\ pointwise samples.  This is addressed in detail in the companion work \cite{BADHFramesPart2}.  
Second, we have used a fixed parameter $\epsilon$ throughout.  However, adaptive (i.e.\ function dependent) choices would likely be beneficial.  This warrants further study.  Third, since the focus of this paper has been on general frames, we have not considered fast computations (which is more specific to the frame employed).  A fast algorithm for computing Fourier extensions in the one-dimensional setting was introduced in \cite{LyonFast} for the special case where the extension interval has exactly twice the length of the original interval.  A more recent alternative that also generalizes to higher dimensions is described in \cite{matthysen2015fastfe,matthysen2017fastfe2d}.  Crucial elements in the latter approach include the link to the theory of bandlimited functions, the special prolate spheroidal wave functions and a phenomenon called the \emph{plunge region} in sampling theory. It seems these elements may generalize to other types of frames, which is a topic that will be considered in future work.  Fourth and finally, we have not discussed the accuracy of computing the SVD of $\bm{G}_N$ in floating-point arithmetic, and its effect on the numerical projection.  Numerical experiments suggest this does not have a substantial impact on the error, but since $\bm{G}_N$ is severely ill-conditioned a careful analysis should be carried out.  We expect the structure of the singular values (in particular, their tendency to divide into `good' singular value away from zero and `bad' singular values near zero) is important in this regard.

\section*{Acknowledgements}
The initial ideas for this paper were first discussed during the Research Cluster on ``Computational Challenges in Sparse and Redundant Representations'' at ICERM in November 2014.  The authors would like to thank all the participants for the useful discussions and feedback received during the program.  They would also like to thank John Benedetto, Juan Manuel C\'ardenas, Pete Casazza, Vincent Copp{\'e}, Roel Matthysen, Sebasti\'an Scheuermann, Mikael Slevinsky, Thomas Strohmer, Nick Trefethen, Andy Wathen and Marcus Webb.  The first author is supported by NSERC grant 611675, as well as an Alfred P.\ Sloan Research Fellowship. The second author is supported by FWO-Flanders projects G.0641.11 and G.A004.14, as well as by KU Leuven project C14/15/055.

\bibliographystyle{abbrv}
\small
\renewcommand{\baselinestretch}{1.25}
\bibliography{FramesStabilityRefs}

\end{document}